# PROOF OF A COMBINATORIAL CONJECTURE POSED IN "THE BLIMPY SHAPE OF HEADY-S AND TAILY-S BIT STRINGS"

**Bruce Levin** BL6@COLUMBIA.EDU
*Department of Biostatistics*
*Mailman School of Public Health*
*Columbia University*
*New York, NY 10032, USA*

**Abstract**

We demonstrate three properties conjectured to hold for a certain function by Levin (2025) in a study of the "blimpy" graphical shape of the number of bit strings with a given score under an interesting scoring system. The properties include discrete convexity, a simple formula for the greatest argument at which the function is negative, and a positive expectation under a certain probability function. A new set of inequalities which imply the latter is presented and proved under some monotonicity assumptions.

## 1. Introduction: statement of problem, motivation, and synopsis of new results

Let $l \ge 6$ be a given integer and let $\tau = \lfloor l/3 \rfloor \ge 2$. For $j = 1,\ldots,\tau$, define the discrete functions

$$g_l(j) := \frac{(l+j)^{(2j)}(l+j)^{(2j)}}{(l+3j)^{(2j)}(l-j)^{(2j)}}, \tag{1.1}$$

$$h_l(j) := \frac{l+j}{l-j}, \tag{1.2}$$

and

$$d_l(j) := g_l(j) - h_l(j), \tag{1.3}$$

where the exponents in (1.1) denote descending factorial powers, $x^{(m)} = x(x-1)\cdots(x-m+1)$ for integer $m$ with $x^{(0)} = 1$. These functions arose in a study by Levin (2025, hereinafter "L25") of the "blimpy" shape of the graph of the logarithm of the number of bit strings of length $n$ which possess certain "score" values, when plotted against those values, under an interesting scoring system for bit strings. In Property 6.3 of that article, we stated without proof that for $l \ge 6$, $d_l(j)$ has positive second-order finite differences; gave a simple formula for where the last negative value of $d_l(j)$ occurs, namely at $a_l := \max\{j : d(j) < 0\} = \lfloor \sqrt{l}/2 \rfloor$; and, most importantly, asserted that $a_l$ satisfies a certain expectation inequality, to wit,

$$\frac{E_l[d(J)\,|\,J > a_l]}{E_l[-d(J)\,|\,J \le a_l]} > \frac{P_l[J \le a_l]}{P_l[J > a_l]}, \tag{1.4}$$

where $J$ is a discrete random variable taking values $j = 1,...,\tau$ with probability function $Q_l(j) = P_l[J = j]$ defined below at (1.10) and where the expectations are taken with respect to $Q_l(j)$. By cross-multiplication and the rule of total expectation, inequality (1.4) means that the unconditional expectation of $d(J)$ is positive, which was a key objective in L25 for reasons explained below. The purpose of the present report is to establish the convexity property with Theorem 1.1 below; to prove the validity of the formula for $a_l$ with Theorem 1.2; and to demonstrate (1.4) in the most stringent case under certain monotonicity assumptions with Theorem 1.3.

We note that functions (1.1) to (1.3) are special cases of more general functions defined at equations (6.10)-(6.12), respectively, in L25. There, we had used the subscript $n$ to index the more general functions and in Lemma 6.11, we derived formula (1.1) from (6.10) in the special case $n=4l$. Here we use the subscript $l$ instead of $n$ since we focus solely on that special case. We had also set $h_n(0) = -1$ in order to be consistent at $j=0$ with the more general definition (6.11), namely, $h_n(j) = \{(l+j)/(l-j)\}\{(j-r/4)/(j+r/4)\}$, when $n = 4l + r$ with $r=1$, 2, or 3. Here, though, we'll use the more natural value $h_l(0) = 1$ with no contradiction or loss of generality because any assertion we make below about the convexity of $d_l(j)$ as a function of $j$ using the value $h_l(0) = 1$ would hold *a fortiori* if the value –1 were used instead. We briefly return to the case $r > 0$ in the discussion Section 5.

The motivation for studying the properties of $d_l(j)$ begins with the following probability puzzle posed by Daniel Litt of the University of Toronto.

> Alice and Bob flip a (fair) coin 100 times. Anytime there are two heads in a row, Alice gets a point; when a head is followed by a tail, Bob gets a point. So in the sequence THHHT, Alice gets two points and Bob gets one. Who is more likely to win?
>
> [Levin (2024), paraphrasing Klarreich (2024)]

In an effort to clarify the counterintuitive answer—that Bob is more likely to win than Alice—we undertook a study of

> …the "blimpy" shape of the graph of the logarithm of the number of bit strings of length *n* which possess certain "score" values, when plotted against those values, under an interesting scoring system for bit strings.
>
> [L25]

We call a binary string $x^{(n)} = (x_1,...,x_n) \in \{0,1\}^n$ with $x_n = 1$ a *heady* bit string (of length *n*) or a *taily* bit string if $x_n = 0$. The *score* in question is the net difference between Alice's and Bob's point totals after $n$ tosses, namely, $S(x^{(n)}) = \sum_{i=1}^{n-1} I[x_i = x_{i+1} = 1] - \sum_{i=1}^{n-1} I[x_i = 1,\ x_{i+1} = 0]$. For a proof of why Bob has a winning advantage, see Levin (2024); Section 1 of L25 further discusses the counter-intuitive result.

For a given integer $s$, we call a heady (respectively, taily) bit string with score $S(x^{(n)}) = s$ a *heady-s* (respectively, *taily-s*) bit string. Let $H_s(n)$ denote the number of heady-$s$ bit strings of length $n$. Levin (2024) and L25 give explicit formulas for $H_s(n)$. In particular, when $s=0$ or $s=-1$, we have

$$H_0(n) = 1 + \sum_{k=1}^{\lfloor (n-1)/3 \rfloor} \binom{2k}{k}\binom{n-1-2k}{k}$$

and

$$H_{-1}(n) = \sum_{k=1}^{\lfloor n/3 \rfloor} \binom{2k-1}{k}\binom{n-2k}{k}.$$

A key goal of L25 was to demonstrate that $s=0$ is a modal value of $H_s(n)$ and part of that requires proving that $H_0(n) \ge H_{-1}(n)$ for all $n \ge 1$. Lemma 2.5 and Remark 2.1 of L25 provided a necessary and sufficient condition for the slightly stronger result that when $n \ge 7$, $H_0(n) > H_{-1}(n)$. The condition, which has some combinatorial interest in its own right, can be stated as follows. For given $n \ge 7$, let $\kappa = \lfloor n/3 \rfloor$ and let $K$ be a discrete random variable over integers $k \in \{1,...,\kappa\}$ with probabilities proportional to $\binom{2k-1}{k}$, i.e., $P_\kappa[K=k] = \binom{2k-1}{k} \Big/ \sum_{i=1}^{\kappa}\binom{2i-1}{i}$. *Then* $H_0(n) > H_{-1}(n)$ *if and only if* $2E_\kappa[\binom{n-1-2K}{K}] \ge E_\kappa[\binom{n-2K}{K}]$ or equivalently by Pascal's Triangle identity,

$$E_\kappa\left[\binom{n-1-2K}{K}\right] \ge E_\kappa\left[\binom{n-1-2K}{K-1}\right], \tag{1.5}$$

where the expectations are taken with respect to $P_\kappa[K=k]$. [As an aside, strict inequality also holds for $n=2$, 4, and 5, but $H_0(n) = H_{-1}(n)$ for $n=1$, 3, and 6. In any case, (1.5) implies $H_0(n) \ge H_{-1}(n)$ for all $n \ge 1$.] *So our motivation for the present work was to demonstrate (1.5), which we will do in the special case when n is a multiple of 4 under some remarkable monotonicity properties.*

There is a logical connection between inequalities (1.5) and (1.4) which we now sketch; see L25 for the details. For the time being, we keep $n \ge 7$ arbitrary until equation (1.10) below. First, we apply a change of measure from $P_\kappa[K=k]$ to probability function (p.f.)

$$P_\kappa^*[K=k] \propto \binom{2k-1}{k-1}\binom{n-1-2k}{k-1} \quad \text{for } k = 1,...,\kappa = \lfloor n/3 \rfloor \ge 2 \text{ for } n \ge 7. \tag{1.6}$$

Then defining $l_n = n/4$, it is easy to show that (1.5) is equivalent to

$$E_\kappa^*[K^{-1}] \ge l_n^{-1} \tag{1.7}$$

where the expectation is with respect to $P_\kappa^*$. Now direct calculation verifies that (1.7) holds for each $n$=7,…,23, in fact with *strict* inequality except for $n$=7. So, without loss of generality we assume that $n \ge 24$ with $l$, the integer part of $l_n$, satisfying $l = \lfloor n/4 \rfloor \ge 6$, and take up the stronger goal of showing

$$E_\kappa^*[K^{-1}] > l_n^{-1}. \tag{1.8}$$

Lemma 6.9 of L25 then shows that for $n \ge 24$, (1.8) already holds if it holds *conditionally* on a certain symmetric subset centered on $l$, to wit,

$$E_\kappa^*[K^{-1} \mid C_n] > l_n^{-1} \quad \text{where} \quad C_n = [2l - \kappa \le k \le \kappa]. \tag{1.9}$$

This allows us to focus on a new random variable $J \in \{0,\ldots,\tau\}$, where $\tau = \kappa - l = \lfloor n/3 \rfloor - \lfloor n/4 \rfloor \ge 2$, with p.f. proportional to

$$\frac{(j + r/4)(l - j)}{(l + j)^2} P_\kappa^*[K = l + j] \quad \text{for } j = 0,\ldots,\tau,$$

where $r = n \bmod 4$ with $n = 4l + r$.

Now restricting attention to the special case $n$=4$l$, $r$=0, we write the above p.f. as

$$Q_l(j) = P_l[J = j] \propto \frac{j(l-j)}{(l+j)^2} P_\kappa^*[K = l + j] \propto \frac{j(l-j)}{(l+j)^2} \binom{2l + 2j - 1}{l + j - 1} \binom{2l - 2j - 1}{l + j - 1} \tag{1.10}$$

for $j = 1,\ldots,\tau = \lfloor l/3 \rfloor$, which we may do as $j = 0$ has zero probability under $Q_l(j)$. *Then Lemmas 6.10 and 6.11 of L25 show that for $n \ge 24$ with $\tau \ge 2$, (1.9) holds if and only if*

$$E_l[d_l(J)] > 0, \tag{1.11}$$

where the expectation is taken with respect to $Q_l$. Lemma 6.12 of L25 further shows that

$$d_l(1) < 0 < d_l(\tau). \tag{1.12}$$

In order to show (1.11) holds, as mentioned at the start, L25 stated without proof that $d_l(j)$ satisfies the following.

<u>Property (6.3) of L25</u>. For any given integer $l \ge 6$ and $j = 1,\ldots,\tau = \lfloor l/3 \rfloor \ge 2$, the following three properties hold.

(i) [convexity]. $d_l(j)$ *has increasing first-order differences*, which we write here as

$$\Delta d_l(j) = d_l(j+1) - d_l(j). \tag{1.13}$$

(ii) [formula for $a_l$]. From (1.12) and (1.13), $d_l(j)$ has a maximum index for which it is negative, call it $a_l = \max\{j : d_l(j) < 0\}$, with $1 \le a_l < \tau$. *Then* $a_l$ *has the simple formula*

$$a_l = \lfloor \sqrt{l}/2 \rfloor. \tag{1.14}$$

(iii) [positive expectation for $d_l(J)$].

$$\frac{E_l[d_l(J) \mid J > a_l]}{E_l[-d_l(J) \mid J \le a_l]} > \frac{P_l[J \le a_l]}{P_l[J > a_l]}. \tag{1.15}$$

As mentioned after (1.4), inequality (1.15) [the same as (1.4)] implies (1.11) by cross-multiplication and the rule of total expectation. *Thus, by establishing Property (6.3) of L25, we conclude that sufficient conditions (1.9) and (1.8) hold, whence we obtain both inequality (1.5) and* $H_0(n) > H_{-1}(n)$.

For the remainder of this section, we synopsize the new results which establish Property 6.3 of L25.

Theorem 1.1. The following assertions hold for each integer $l \ge 6$.

*(i)* $h_l(j)$ *increases in j and has strictly positive second-order finite differences,*

$$\Delta^2 h_l(j) = h_l(j+1) - 2h_l(j) + h_l(j-1) > 0 \quad \textit{for} \quad j = 1,\ldots,\tau-1. \tag{1.16}$$

*(ii)* $\Delta^2 h_l(j)$ *increases in j as well, i.e.,*

$$\Delta^3 h_l(j) = \Delta^2 h_l(j+1) - \Delta^2 h_l(j) > 0 \quad \textit{for} \quad j = 1,\ldots,\tau-2. \tag{1.17}$$

*(iii)* $g_l(j)$ *increases in j and also has strictly positive second-order finite differences,*

$$\Delta^2 g_l(j) = g_l(j+1) - 2g_l(j) + g_l(j-1) > 0 \quad \textit{for} \quad j = 1,\ldots,\tau-1. \tag{1.18}$$

*(iv)* $\Delta^2 g_l(j)$ *increases in j as well, i.e.,*

$$\Delta^3 g_l(j) = \Delta^2 g_l(j+1) - \Delta^2 g_l(j) > 0 \quad \textit{for} \quad j = 1,\ldots,\tau-2. \tag{1.19}$$

*(v) The smallest value of* $\Delta^2 g_l(j)$ *exceeds the largest value of* $\Delta^2 h_l(j)$, *i.e.,*

$$\Delta^2 g_l(1) > \Delta^2 h_l(\tau-1). \tag{1.20}$$

Theorem 1.1 implies that

$$\Delta^2 d_l(j) = d_l(j+1) - 2d_l(j) + d_l(j-1) > 0 \quad \text{for} \quad j = 1,\ldots,\tau-1, \tag{1.21}$$

for by (iv), (v), and (ii), $\Delta^2 g_l(j) \ge \Delta^2 g_l(1) > \Delta^2 h_l(\tau-1) \ge \Delta^2 h_l(j)$ and since $\Delta^2$ is a linear functional, $\Delta^2 d_l(j) = \Delta^2\{g_l(j) - h_l(j)\} = \Delta^2 g_l(j) - \Delta^2 h_l(j) > 0$. Inequality (1.21) means $d_l(j)$ has increasing first-order finite differences as asserted in Property 6.3 of L25.

Before continuing, we should clarify a matter of notation. We use the usual definition of first *forward* differences, e.g., $\Delta h_l(j) := h_l(j+1) - h_l(j)$, but we have chosen to use *j* rather than *j*–1 in the notation for second-order finite differences *centered* at *j*, in order to allow the starting index to be fixed at $j = 1$. For example, we write $\Delta^2 h_l(j) = \Delta h_l(j) - \Delta h_l(j-1) = \Delta\{\Delta h_l(j-1)\}$ for $j = 1,...,\tau-1$ rather than $\Delta\{\Delta h_l(j)\}$ for $j = 0,...,\tau-2$. For third-order differences, we still have the identity $\Delta^3 h_l(j) = \Delta\{\Delta^2 h_l(j)\} = \Delta^2\{\Delta h_l(j)\}$ for $j = 1,...,\tau-2$, because

$$\Delta\{\Delta^2 h_l(j)\} = \Delta^2 h_l(j+1) - \Delta^2 h_l(j) = \{h_l(j+2) - 2h_l(j+1) + h_l(j)\} - \{h_l(j+1) - 2h_l(j) + h_l(j-1)\}$$

$$= \Delta h_l(j+1) - 2\Delta h_l(j) + \Delta h_l(j-1) = \Delta^2\{\Delta h_l(j)\}.$$

The next result shows that $a_l$ is given by the integer part of $\sqrt{l}/2$.

<u>Theorem 1.2</u>. *(i) For values of* $l = 4a^2$ $(a = 1, 2,...)$, *asymptotically as* $a \to \infty$,

$$\log g_{4a^2}(a) = \frac{1}{2a} - \frac{3}{32a^3} + O(a^{-5}) \tag{1.22}$$

*and*

$$\log h_{4a^2}(a) = \frac{1}{2a} + \frac{1}{96a^3} + O(a^{-5}). \tag{1.23}$$

*Thus ignoring terms of order* $a^{-5}$, $g_{4a^2}(a) < h_{4a^2}(a)$.

*(ii) For the same values of l, but evaluating j at* $a+1$ *rather than a, as* $a \to \infty$,

$$\log g_{4a^2}(a+1) = \frac{1}{2a} + \frac{3}{2a^2} + \frac{47}{32a^3} + \frac{19}{32a^4} + O(a^{-5}) \tag{1.24}$$

*and*

$$\log h_{4a^2}(a+1) = \frac{1}{2a} + \frac{1}{2a^2} + \frac{1}{96a^3} + \frac{1}{32a^4} + O(a^{-5}). \tag{1.25}$$

*Thus ignoring terms of order* $a^{-5}$, $g_{4a^2}(a+1) > h_{4a^2}(a+1)$.

*(iii) For values of* $l = 4a^2 - 1$ $(a = 1, 2,...)$, *asymptotically as* $a \to \infty$,

$$\log g_{4a^2-1}(a) = \frac{1}{2a} + \frac{7}{32a^3} + O(a^{-5}) \tag{1.26}$$

*and*

$$\log h_{4a^2-1}(a) = \frac{1}{2a} + \frac{13}{96a^3} + O(a^{-5}). \tag{1.27}$$

*Thus ignoring terms of order* $a^{-5}$, $g_{4a^2-1}(a) > h_{4a^2-1}(a)$.

*(iv) For any fixed* $a \geq 3$, *both* $g_l(a)/h_l(a)$ *and* $g_l(a+1)/h_l(a+1)$ *decrease strictly in l for l in the interval* $4a^2 \leq l \leq 4(a+1)^2 - 1$.

Assertions (i) and (ii) of Theorem 1.2 establish that, ignoring terms of order $a^{-5}$ and higher, $a_l = \lfloor \sqrt{l}/2 \rfloor$ holds in the base case $l = 4a^2$. Assertion (iv) then shows that when $a \ge 3$, $g_l(a)/h_l(a) < g_{4a^2}(a)/h_{4a^2}(a) < 1$ for all other $l$ in the interval $4a^2 \le l \le 4(a+1)^2 - 1$, whence $g_l(a) < h_l(a)$, and over the same interval, $g_l(a+1)/h_l(a+1) > g_{4(a+1)^2-1}(a+1)/h_{4(a+1)^2-1}(a+1) > 1$ by (iii) (applied to $a+1$), whence $g_l(a+1) > h_l(a+1)$. Thus $a_l = \lfloor \sqrt{l}/2 \rfloor$ holds for all $l \ge 4 \cdot 3^2 = 36$ for sufficiently large $a$. Direct calculation then confirms that the three inequalities $g_{4a^2}(a) < h_{4a^2}(a)$, $g_{4a^2}(a+1) > h_{4a^2}(a+1)$, and $g_{4a^2-1}(a) < h_{4a^2-1}(a)$ hold for all $a$=3,...,100. Calculation also shows $g_l(1) < h_l(1)$ while $g_l(2) > h_l(2)$ for $l = 6,...,15$ and $g_l(2) < h_l(2)$ while $g_l(3) > h_l(3)$ for $l = 16,...,35$.

The next result implies the key inequality (1.4) in the special case $n$=4$l$ under consideration. When $l$=6, 7, or 8 with $\tau = 2$, (1.4) reduces to $d_l(2)/\{-d_l(1)\} > Q_l(1)/Q_l(2)$, which is easily confirmed by direct calculation. So without loss of generality we may assume $l \ge 9$ with $\tau \ge 3$, in which case we have $2a_l + 1 \le \tau$. This is because $a_l = 1$ for $9 \le l \le 15$ so $2a_l + 1 = 3 \le \tau$, while for $l$=16 or 17, $a$=2 and $\tau = 5$, so $2a_l + 1 = 5 = \tau$. For $l \ge 18$ and $\tau \ge 6$, write $l = 3\tau + \delta$ with $\delta \in \{0, 1, 2\}$. Then $2a_l + 1 = 2\lfloor \sqrt{l}/2 \rfloor + 1 = 2\lfloor \sqrt{3\tau + \delta}/2 \rfloor + 1 = 2\lfloor \sqrt{(3/4)\tau + \delta/4} \rfloor + 1$ which is no greater than $\tau$ if $\sqrt{(3/4)\tau + \delta/4} \le (\tau - 1)/2$, which holds if and only if $3\tau + \delta \le \tau^2 - 2\tau + 1$, i.e., $\tau^2 - 5\tau + 1 - \delta \ge 0$, which holds for $\tau \ge 6$.

Now let $J$ take values $j = 1,...,\tau$ with probabilities given in (1.10), namely,

$$Q_l(j) := \frac{j(l-j)}{(l+j)^2}\binom{2(l+j)-1}{l+j-1}\binom{4l-1-2(l+j)}{l+j-1} \Big/ \sum_{i=1}^{\tau} \frac{i(l-i)}{(l+i)^2}\binom{2(l+i)-1}{l+i-1}\binom{4l-1-2(l+i)}{l+i-1}. \tag{1.28}$$

<u>Theorem 1.3</u>. *For* $l \ge 9$ *and each* $j = 1,...,a_l = \lfloor \sqrt{l}/2 \rfloor$,

$$d_l(a_l + 1 + j)Q_l(a_l + 1 + j) \ > \ -d_l(j)Q_l(j). \tag{1.29}$$

We call (1.29) the *pairwise inequalities for* $d_l(l)Q_l(l)$. They imply (1.4) and thus $E_l[d_l(J)] > 0$, since

$$E_l[d_l(J) \mid J > a]P_l[J > a] = \sum_{j=a_l+1}^{\tau} d_l(j)Q_l(j) > \sum_{j=a_l+2}^{2a_l+1} d_l(j)Q_l(j) > -\sum_{j=1}^{a_l} d_l(j)Q_l(j) = -E_l[d_l(J) \mid J \le a]P_l[J \le a],$$

where the first inequality holds because both $d_l(j) > 0$ and $Q_l(j) > 0$ for $j = a_l + 1$ and $2a_l + 1 < j \le \tau$ in cases where $2a_l + 1 < \tau$.

We prove Theorem 1.1 in the next section, Theorem 1.2 in Section 3, and Theorem 1.3 in Section 4.

## 2. Proof of Theorem 1.1.

Proof of assertion 1.1(i). $h_l(j) = (l+j)/(l-j)$ increases in $j$ and $\Delta^2 h_l(j) > 0$ for $j = 1,\ldots,\tau-1$ because the function $h_l(x) = (l+x)/(l-x)$ for continuous $x$ is strictly convex in $x$ for $x \in [0,\tau]$. The positivity of $\Delta^2 h_l(j)$ follows from Jensen's inequality, for if $J$ is a random variable with $P[J = j+1] = P[J = j-1] = ½$, then $½\Delta^2 h_l(j) = E[h_l(J)] - h_l(E[J]) > 0$. Thus assertion 1.1(i) holds.

Proof of assertion 1.1(ii). We first simplify $\Delta^2 h_l(j) = \dfrac{l+j+1}{l-j-1} - 2\left(\dfrac{l+j}{l-j}\right) + \dfrac{l+j-1}{l-j+1}$ with the result in (2.1) below.

$$
\begin{aligned}
(l-j+1)^{(3)}\Delta^2 h_l(j) &= (l+j+1)(l-j+1)(l-j) - 2(l+j)(l-j+1)(l-j-1) + (l+j-1)(l-j)(l-j-1) \\
&= (j+l+1)\{j-(l+1)\}(j-l) \;-\; 2(j+l)\{j-(l+1)\}\{j-(l-1)\} \;+\; (j+l-1)(j-l)\{j-(l-1)\} \\
&= \{j^2-(l+1)^2\}(j-l) \;-\; 2(j+l)\{j^2-2lj+l^2-1\} \;+\; (j-l)\{j^2-(l-1)^2\} \\
&= \{j^3-lj^2-(l+1)^2 j+l(l+1)^2\} \;-\; 2\{j^3-lj^2-(l^2+1)j+l(l^2-1)\} \;+\; \{j^3-lj^2-(l-1)^2 j+l(l-1)^2\}.
\end{aligned}
$$

The coefficients multiplying $j^3$ and $j^2$ happily vanish, and since $-(l+1)^2 + 2(l^2+1) - (l-1)^2 = 0$, so does $coeff(j)$. The constant term is $l(l+1)^2 - 2l(l^2-1) + l(l-1)^2 = l\{(l+1)-(l-1)\}^2 = 4l$. Therefore,

$$
\Delta^2 h_l(j) = \frac{4l}{(l-j+1)^{(3)}}, \tag{2.1}
$$

which clearly increases in $j$. In particular, for $j = 1,\ldots,\tau-2$,

$$
\Delta^3 h_l(j) = \Delta^2 h_l(j+1) - \Delta^2 h_l(j) = \frac{4l}{(l-j)(l-j-1)}\left\{\frac{1}{l-j-2} - \frac{1}{l-j+1}\right\} = \frac{12l}{(l-j+1)^{(4)}} > 0 .
$$

This establishes assertion 1.1(ii).

We note for use in assertion 1.1(v) below that, from (2.1), the largest value of $\Delta^2 h_l(j)$ is

$$
\Delta^2 h_l(\tau-1) = \frac{4l}{(l-\tau+2)^{(3)}} \text{ with } \tau = \lfloor l/3 \rfloor \le l/3 . \tag{2.2}
$$

Also, since $l-\tau \ge 2l/3$, we have $(l-\tau+2)^{(3)} \ge \left(\dfrac{2l}{3}+2\right)\left(\dfrac{2l}{3}+1\right)\left(\dfrac{2l}{3}\right) = \left(\dfrac{8l^3}{27}\right)\left(1+3l^{-1}\right)\left(1+\dfrac{3}{2}l^{-1}\right)$,

so an upper bound for $\Delta^2 h_l(\tau-1)$ as a function of $l$ is

$$\Delta^2 h_l(\tau-1) \le \frac{(27/2)l^{-2}}{(1+3l^{-1})\{1+(3/2)l^{-1}\}} < (27/2)l^{-2} \tag{2.3}$$

with exact equality on the left if $l$ is a multiple of 3.

<u>Proof of assertion 1.1(iii)</u>. It will be convenient to re-express $g_l(j)$ for positive $j$ as follows.

$$g_l(j) = \frac{(l+j)^{(2j)}(l+j)^{(2j)}}{(l+3j)^{(2j)}(l-j)^{(2j)}} = \prod_{i=0}^{2j-1} \frac{(l+j-i)^2}{(l+j-i+2j)(l+j-i-2j)^2} = \prod_{i=0}^{2j-1} \frac{(l+j-i)^2}{(l+j-i)^2-4j^2}$$

$$= \prod_{i=0}^{2j-1}\left\{1+\frac{4j^2}{(l+j-i)^2-4j^2}\right\} = \prod_{i=0}^{2j-1}\left\{1+\frac{1}{(\frac{l+j-i}{2j})^2-1}\right\} = \prod_{i=0}^{2j-1}\left[1+\frac{1}{\left\{\frac{1+(l-i)/j}{2}\right\}^2-1}\right]. \tag{2.4}$$

For given constants $0 \le a < 1$ and $b > 0$, define the real function

$$f_{a,b}(x) := \frac{1}{(a+b/x)^2-1} \tag{2.5}$$

for $0 < x < b/(1-a)$, so that $f_{a,b}(x) > 0$. Then $g_l(j)$ in (2.4) is of the form $g_l(j) = \prod_{i=0}^{2j-1}\{1+f_{a,b_i}(j)\}$ with $a=½$ and $b_i = (l-i)/2$. In this case we use the abbreviation $f_i(x) = f_{1/2,\,(l-i)/2}(x)$ for $0 < x < b_i/(1-a) = l-i$ and write

$$g_l(j) = \prod_{i=0}^{2j-1}\left\{1+f_i(j)\right\}. \tag{2.6}$$

Now $f_i(x)$ increases in $x$ for each $i$, and this implies $g_l(j)$ increases in $j$. That is because

$$g_l(j+1) = \prod_{i=0}^{2j+1}\{1+f_i(j+1)\} = \left\{1+f_{2j+1}(j+1)\right\}\left\{1+f_{2j}(j+1)\right\}\prod_{i=0}^{2j-1}\{1+f_i(j+1)\}$$

$$> \prod_{i=0}^{2j-1}\left\{1+f_i(j+1)\right\} > \prod_{i=0}^{2j-1}\left\{1+f_i(j)\right\} = g_l(j).$$

[As a technical detail, we confirm that for any $j \le \tau-1$ and any $i \le 2j+1$, we have $j+1 < l-i$, so that $j+1$ is in the domain of $f_i(x)$, with $f_i(j+1)$ defined and positive. This holds because $3j \le 3\tau - 3 = 3\lfloor l/3 \rfloor - 3 \le l-3 < l-2$, whence $j+1 < l-2j-1 \le l-i$.]

Next we show $\Delta^2 g_l(j) > 0$ for $j = 1,\ldots,\tau-1$. This might seem apparent insofar as $f_i(x)$ is convex in $x$ for each $i$ (as shown below), and the product of positive, increasing, convex functions is convex. However, $j$ appears in (2.6) in both the product's factors and its upper limit, and though this turns out not to be an issue, we will account for it to be rigorous.

The following notation will be helpful. For given integers $n$ and $k$ with $n \geq 3$ and $0 \leq k \leq n$, let $I_{nk} = \{(i) = (i_1,...,i_k) \, : \, 0 \leq i_1 < \cdots < i_k \leq n\}$ be the set of *ordered* distinct-integer $k$-tuples $\in \{0,...,n\}$.

For $(i) \in I_{nk}$, define $f_{(i)}(j) = \prod_{\alpha=1}^{k} f_{i_\alpha}(j)$. Then with $n = 2j-1$, we can write $g_l(j) = \prod_{i=0}^{2j-1} \{1 + f_i(j)\} = 1 + \sum_{(i)\in I_{n1}} f_{(i)}(j) + \cdots + \sum_{(i)\in I_{n,n+1}} f_{(i)}(j)$ where the first sum equals $\sum_{i=0}^{n} f_i(j)$ and the last sum equals $f_0(j) \cdots f_n(j)$.

For the gist of the proof, first consider the case $j=1$ for which we have

$$g_l(2) = \prod_{i=0}^{3} \{1 + f_i(2)\} = 1 + \sum_{i=0}^{3} f_i(2) + \sum_{(i)\in I_{32}} f_{(i)}(2) + \sum_{(i)\in I_{33}} f_{(i)}(2) + \sum_{(i)\in I_{34}} f_{(i)}(2) \, .$$

Thus $g_l(2)$ exceeds $1 + 2\{f_0(2) + f_1(2) + f_0(2) f_1(2)\}$ because (a) $f_i(j) = \dfrac{1}{\{½ + ½(l-i)/j\}^2 - 1}$ increases in $i$ for each fixed $j$, so that $f_2(2) + f_3(2) > f_0(2) + f_1(2)$, whence $\sum_{i=0}^{3} f_i(2) > 2\{f_0(2) + f_1(2)\}$;

(b) each of the six products $f_{(i)}(2) = f_{i_1}(2) f_{i_2}(2)$ for $(i) \in I_{32}$ is no less than $f_0(2) f_1(2)$; and (c) the remaining sums for $(i) \in I_{33}$ and $(i) \in I_{34}$ are positive. But now

$$g_l(1) = \prod_{i=0}^{1} \{1 + f_i(1)\} = 1 + f_0(1) + f_1(1) + f_0(1) f_1(1) \, ,$$

and with $g_l(0) = 1$, we find

$$\Delta^2 g_l(1) = g_l(2) - 2 g_l(1) + g_l(0) > 2\{f_0(2) + f_1(2) + f_0(2) f_1(2)\} \; - \; 2\{f_0(1) + f_1(1) + f_0(1) f_1(1)\} \; > \; 0$$

since $f_i(j)$ also increases in $j$ for each fixed $i$.

The proof of $\Delta^2 g_l(j) > 0$ for $j > 1$ is similar. We break $\Delta^2 g_l(j)$ into *initial sums*, involving only $(i) \in I_{pk}$ with $p = 2(j-1) - 1 = 2j - 3$, plus *additional sums*, involving either $(i) \in I_{mk} \setminus I_{pk}$ with $m = 2j+1$, arising from $g_l(j+1)$, or $(i) \in I_{nk} \setminus I_{pk}$ with $n = 2j-1$, arising from $g_l(j)$. We first show that the contribution of the additional sums to $\Delta^2 g_l(j)$ is positive, i.e., the sum of all terms involving one or more subscripts $i_\alpha > 2j-3$ in $\Delta^2 g_l(j)$ only augments the initial sums in $\Delta^2 g_l(j)$, which involve only subscripts $i_\alpha \leq 2j-3$. We then show that those initial sums in $\Delta^2 g_l(j)$ are positive. We have

$$g_l(j+1)=\prod_{i=0}^{2j+1}\{1+f_i(j+1)\}=1+\sum_{i=0}^{m}f_i(j+1)+\cdots+\sum_{(i)\in I_{mk}}f_{(i)}(j+1)+\cdots$$

$$=1+\sum_{i=0}^{p}f_i(j+1)+\cdots+\sum_{(i)\in I_{pk}}f_{(i)}(j+1)+\cdots \qquad \text{(initial sums)}$$

$$\sum_{i=2j-2}^{2j+1}f_i(j+1)+\cdots+\sum_{(i)\in I_{mk}\setminus I_{pk}}f_{(i)}(j+1)+\cdots \quad \text{(additional sums).}$$

Similarly,

$$g_l(j)=\prod_{i=0}^{2j-1}\{1+f_i(j)\}=1+\sum_{i=0}^{n}f_i(j)+\cdots+\sum_{(i)\in I_{nk}}f_{(i)}(j)+\cdots$$

$$=1+\sum_{i=0}^{p}f_i(j)+\cdots+\sum_{(i)\in I_{pk}}f_{(i)}(j)+\cdots \qquad \text{(initial sums)}$$

$$+\sum_{i=2j-2}^{2j-1}f_i(j)+\cdots+\sum_{(i)\in I_{nk}\setminus I_{pk}}f_{(i)}(j)+\cdots \quad \text{(additional sums).}$$

In $g_l(j+1)$, for each $k=1,\ldots,m+1$ and $r=1,\ldots,4\wedge k$, suppose $(i)\in I_{mk}\setminus I_{pk}$ has $r$ components greater than 2$j$–3 and $k-r\geq 0$ components no greater than 2$j$–3. Then for each fixed $(i')\in I_{p,k-r}$, there are $\binom{4}{r}$ values of $(i)$ concatenating $(i')$ and $(i_{k-r+1},\cdots,i_k)$ with $2j-2\leq i_{k-r+1}<\cdots<i_k\leq 2j+1$. In particular, there are four such $(i)$ when $r=1$ and six when $r=2$ (and $k\geq 2$). When $r=1$, $f_{(i)}(j+1)=\left\{\prod_{\alpha=1}^{k-1}f_{i'_\alpha}(j+1)\right\}f_{i_k}(j+1)$, and with

$$f_{2j-2}(j+1)+\cdots+f_{2j+1}(j+1)\ >\ 2\{f_{2j-2}(j+1)+f_{2j-1}(j+1)\}\ >\ 2\{f_{2j-2}(j)+f_{2j-1}(j)\},$$

the sum of $f_{(i)}(j+1)$ over these four $(i)$ with fixed $(i')$ exceeds

$$2\left\{\prod_{\alpha=1}^{k-1}f_{i'_\alpha}(j+1)\right\}\{f_{2j-2}(j)+f_{2j-1}(j)\}\ >\ 2\left\{\prod_{\alpha=1}^{k-1}f_{i'_\alpha}(j)\right\}\{f_{2j-2}(j)+f_{2j-1}(j)\}.$$

Then summing over all $0\leq i'_1<\cdots<i'_{k-r}\leq 2j-3$, we find the sum of $f_{(i)}(j+1)$ over all $(i)\in I_{mk}\setminus I_{pk}$ with $r=1$ exceeds twice the sum of $f_{(i)}(j)$ over all $(i)\in I_{nk}\setminus I_{pk}$ with $r=1$.

Similarly, when $r=2$ (and $k\geq 2$), $f_{(i)}(j+1)=\left\{\prod_{\alpha=1}^{k-2}f_{i_\alpha}(j+1)\right\}f_{i_{k-1}}(j+1)f_{i_k}(j+1)$, and with

$$\sum_{2j-2\leq i_{k-1}<i_k\leq 2j+1}f_{i_{k-1}}(j+1)f_{i_k}(j+1)\ >\ 6f_{2j-2}(j+1)f_{2j-1}(j+1)\ >\ 2f_{2j-2}(j)f_{2j-1}(j),$$

the sum of $f_{(i)}(j+1)$ over the six $(i)$ with fixed $(i')\in I_{p,k-2}$, certainly exceeds

$2\left\{\prod_{\alpha=1}^{k-2} f_{i'_\alpha}(j)\right\} f_{2j-2}(j) f_{2j-1}(j)$. Summing over all $(i')$, we find the sum of $f_{(i)}(j+1)$ over all $(i) \in I_{mk} \setminus I_{pk}$ with $r$=2 exceeds twice the sum of $f_{(i)}(j)$ over all $(i) \in I_{nk} \setminus I_{pk}$ with $r$=2.

Other sums of $f_{(i)}(j+1)$ with $r$=3 or 4 for $k = 3,...,m+1$ occur in $g_l(j+1)$, but note that in $g_l(j)$, $k \le n+1$ and the additional sums can have only $\binom{2}{r}$ values of $(i) \in I_{nk} \setminus I_{pk}$ for each $(i') \in I_{p,k-r}$, so no sums with $r$=3 or 4 (or $k = 2j$ or $2j+1$) occur to subtract from those in $g_l(j+1)$. Thus for $g_l(j+1)$, summing over all $(i') \in I_{p,k-r}$, $r$, and $k \le m+1$, we find

$$g_l(j+1) \;>\; 1+ \text{ initial sums } + \; 2\left\{\sum_{i=2j-2}^{2j-1} f_i(j) \;+\cdots+ \sum_{(i)\in I_{nk}\setminus I_{pk}} f_{(i)}(j) +\cdots+ \sum_{(i)\in I_{n,n+1}\setminus I_{p,n+1}} f_{(i)}(j)\right\}$$

where the additional sums in the last expression are precisely those of $2g_l(j)$. As stated above, therefore, the additional sums in $\Delta^2 g_l(j)$ only augment the initial sums, leading to the inequality

$$\Delta^2 g_l(j) = g_l(j+1) - 2g_l(j) + g_l(j-1)$$

$$> \left\{\sum_{i=0}^{2j-3} f_i(j+1) + \sum_{(i)\in I_{p2}} f_i(j+1) + \cdots\right\} - 2\left\{\sum_{i=0}^{2j-3} f_i(j) + \sum_{(i)\in I_{p2}} f_{(i)}(j) + \cdots\right\} + \left\{\sum_{i=0}^{2j-3} f_i(j-1) + \sum_{(i)\in I_{p2}} f_{(i)}(j-1) + \cdots\right\}$$

$= \sum_{i=0}^{2j-3} \Delta^2 f_i(j) + \sum_{(i)\in I_{p2}} \Delta^2 f_{(i)}(j) + \cdots + \sum_{(i)\in I_{p,p+1}} \Delta^2 f_{(i)}(j)$, involving only initial sums of second-order finite differences of $f_i$ and their products.

It remains to show that $\Delta^2 f_{(i)}(j) > 0$ for each $(i) \in I_{pk}$ and $k \le p+1$. In Appendix A we prove that $f_{a,b}(x)$ in (2.5) is strictly convex in $x$ for $0 \le a < 1$ and $b > 0$ for $x \in (0, b/(1-a))$. That suffices, for then, with random variable $J$ as above with $P[J = j+1] = P[J = j-1] = ½$, we have $½\Delta^2 f_i(j) =$ $½\Delta^2 f_i(j) = E[f_i(J)] - f_i(E[J]) > 0$. Moreover, because products of positive, increasing, convex functions are themselves positive, increasing, and convex, we have for each $(i) \in I_{pk}$,

$$\Delta^2 f_{(i)}(j) = \prod_{\alpha=1}^{k} f_{i_\alpha}(j+1) \;-\; 2\prod_{\alpha=1}^{k} f_{i_\alpha}(j) \;+\; \prod_{\alpha=1}^{k} f_{i_\alpha}(j-1) \;>\; 0\,.$$

Summing over all $(i) \in I_{pk}$ and $k = 1,..., p+1$ then yields $\Delta^2 g_l(j) > 0$. Thus assertion 1.1(iii) holds.

<u>Proof of assertion 1.1(iv)</u>. To show $\Delta^2 g_l(j)$ increases in *j*, it would be natural to try to proceed as in (iii), arguing that the additional sums in $\Delta^3 g_l(j)$ only augment the initial sums restricted to subscripts $(i) \in I_{pk}$, and then appealing to positive third derivatives of $f_i(x)$ to conclude $\Delta^3 g_l(j) > 0$. Unfortunately, $f_i(x)$ does *not* (!) have positive third derivatives for all $x \in (0, l-i)$, so that approach will not work. Instead, we will demonstrate directly that $\Delta g_l(j) = g_l(j+1) - g_l(j)$ has strictly positive second-order finite differences by showing that the ratio $g_l(j+1)/g_l(j)$ is also positive, increasing, and has positive second-order finite differences. This suffices upon writing

$$\Delta g_l(j) = g_l(j+1) - g_l(j) = g_l(j)\{\frac{g_l(j+1)}{g_l(j)} - 1\} \quad \text{for } j = 0, \ldots, \tau - 1$$

and applying the following lemma to the two convex sequences $\{g_l(j) : j = 0, \cdots, \tau - 1\}$ and $\{\frac{g_l(j+1)}{g_l(j)} - 1 : j = 0, \cdots, \tau - 1\}$.

Lemma 2.1. If $\{s_j\}$ and $\{t_j\}$ are two positive sequences that increase ($\Delta s_j = s_{j+1} - s_j > 0$, ditto for $\Delta t_j$) and are convex ($\Delta^2 s_j = \Delta s_j - \Delta s_{j-1} = s_{j+1} - 2s_j + s_{j-1} > 0$, ditto for $\Delta^2 t_j$), then $\Delta\{s_j t_j\} > 0$ and $\Delta^2\{s_j t_j\} > 0$.

Proof. See Appendix B. □

To show $\frac{g_l(j+1)}{g_l(j)} > 0$ increases and has positive second-order finite differences, we write the ratio explicitly as

$$\begin{aligned} \frac{g_l(j+1)}{g_l(j)} &= \frac{(l+j+1)^{(2j+2)}(l+j+1)^{(2j+2)}(l+3j)^{(2j)}(l-j)^{(2j)}}{(l+3j+3)^{(2j+2)}(l-j-1)^{(2j+2)}(l+j)^{(2j)}(l+j)^{(2j)}} \\ &= \frac{(l+j+1)^2 \cdots (l-j)^2 \ \cdot \ (l+3j) \cdots (l+j+1) \ \cdot \ (l-j) \cdots (l-3j+1)}{(l+j)^2 \cdots (l-j+1)^2 \ \cdot \ (l+3j+3) \cdots (l+j+2) \ \cdot \ (l-j-1) \cdots (l-3j-2)} \qquad (2.7) \\ &= \frac{(l+j+1)^2(l-j)^2(l+j+1)(l-j)}{(l+3j+3)^{(3)}(l-3j)^{(3)}} = \frac{(l+j+1)^3(l-j)^3}{(l+3j+3)^{(3)}(l-3j)^{(3)}}. \end{aligned}$$

Now for continuous $x$ satisfying $0 \le x < (l-2)/3$, define the function

$$\begin{aligned} L(x) &:= \log \frac{(l+x+1)^3(l-x)^3}{(l+3x+3)^{(3)}(l-3x)^{(3)}} \\ &= 3\log(x+l+1) + 3\log(l-x) - \log(3x+l+3) - \log(3x+l+2) - \log(3x+l+1) \\ &\quad - \log(l-3x) - \log(l-1-3x) - \log(l-2-3x). \end{aligned}$$

Note that $j \le \tau - 1$ is in the domain of $L(x)$, since $j \le \tau - 1 = \lfloor l/3 \rfloor - 1 = \lfloor (l-3)/3 \rfloor < (l-2)/3$.

Then taking derivatives,

$$(1/3)L'(x) = \frac{1}{x+l+1} - \frac{1}{l-x} - \frac{1}{3x+l+3} - \frac{1}{3x+l+2} - \frac{1}{3x+l+1} + \frac{1}{l-3x} + \frac{1}{l-1-3x} + \frac{1}{l-2-3x}$$

$$= \left(\frac{1}{x+l+1} - \frac{1}{3x+l+1}\right) + \left(\frac{1}{l-3x} - \frac{1}{l-x}\right) + \left(\frac{1}{l-1-3x} - \frac{1}{3x+l+3}\right) + \left(\frac{1}{l-2-3x} - \frac{1}{3x+l+2}\right).$$

But for $x \geq 0$, $x+l+1 \leq 3x+l+1$ and $l-3x \leq l-x$; $l-1-3x < 3x+l+3$ if and only if $6x > -4$, which holds; and ditto for $l-2-3x < 3x+l+2$. Hence $L'(x) > 0$, so that $\frac{g_l(j+2)}{g_l(j+1)} = \exp L(j+1) > \exp L(j) = \frac{g_l(j+1)}{g_l(j)}$. Taking the second derivative,

$$(1/3)L''(x) = \left\{\frac{-1}{(x+l+1)^2} + \frac{3}{(3x+l+1)^2}\right\} + \left\{\frac{3}{(l-3x)^2} - \frac{1}{(l-x)^2}\right\}$$

$$+ \frac{3}{(l-1-3x)^2} + \frac{3}{(3x+l+3)^2} + \frac{3}{(l-2-3x)^2} + \frac{3}{(3x+l+2)^2}.$$

Now, $\frac{3}{(3x+l+1)^2} > \frac{1}{(x+l+1)^2}$ if and only if $x < (l+1)/\sqrt{3}$, which holds for $x < (l-2)/\sqrt{3} < (l+1)/\sqrt{3}$. Similarly, $\frac{3}{(l-3x)^2} > \frac{1}{(l-x)^2}$ if and only if $x < l/\sqrt{3}$, which also holds for $x < (l-2)/\sqrt{3} < l/\sqrt{3}$. Since the remaining terms are all positive, $L(x)$ is convex, hence $\exp L(x)$ is convex, hence $\Delta^2\left\{\frac{g_l(j+1)}{g_l(j)}\right\} > 0$. This establishes assertion 1.1(iv).

<u>Proof of assertion 1.1(v)</u>. We show that for all $l \geq 6$ with $\lambda = l^{-1} < 1/6$, (a) $g_l(1) < 1 + 9\lambda^2$ and (b) $g_l(2) > 1 + (160/3)\lambda^2 > 1 + 53\lambda^2$. Then $\Delta^2 g_l(1) = g_l(2) - 2g_l(1) + 1 > 35\lambda^2$. But from (2.3), $\Delta^2 h_l(\tau - 1) < (27/2)\lambda^2$. Therefore $\Delta^2 g_l(1) > \Delta^2 h_l(\tau - 1)$ for all $l \geq 6$, which is assertion (v).

For (a),

$$g_l(1) = \frac{(l+1)^2 l^2}{(l+3)(l+2)(l-1)(l-2)} = \frac{1+2\lambda+\lambda^2}{(1+3\lambda)(1+2\lambda)(1-\lambda)(1-2\lambda)} = \frac{1+2\lambda+\lambda^2}{(1+2\lambda-3\lambda^2)(1-4\lambda^2)}$$

$$= \frac{1 + \frac{4\lambda^2}{1+2\lambda-3\lambda^2}}{1-4\lambda^2}.$$

But for $\lambda \leq 1/6$, $2\lambda - 3\lambda^2 > 0$ and $\frac{1}{1-4\lambda^2} \leq 9/8$, so $g_l(1) < \frac{1+4\lambda^2}{1-4\lambda^2} = 1 + \frac{8\lambda^2}{1-4\lambda^2} \leq 1 + 9\lambda^2$, which is (a).

For (b), dividing through by $l^6$,

$$g_l(2) = \frac{(l+2)^{(4)}(l+2)^{(4)}}{(l+6)^{(4)}(l-2)^{(4)}} = \frac{(l+2)^2(l+1)^2 l^2 (l-1)^2}{(l+6)(l+5)(l+4)(l+3)(l-2)(l-3)(l-4)(l-5)}$$

$$= \frac{(1+2\lambda)^2(1-\lambda^2)^2}{(1+6\lambda)(1+5\lambda)(1+4\lambda)(1+3\lambda)(1-2\lambda)(1-3\lambda)(1-4\lambda)(1-5\lambda)}$$

$$= \frac{(1+4\lambda+4\lambda^2)(1-\lambda^2)^2}{(1+6\lambda)(1-2\lambda)(1-9\lambda^2)(1-16\lambda^2)(1-25\lambda^2)}.$$

Now $(1-9\lambda^2)(1-16\lambda^2)(1-25\lambda^2) = 1-50\lambda^2+769\lambda^4-3600\lambda^6 < 1-50\lambda^2+769\lambda^4$ and $(1-\lambda^2)^2 > 1-2\lambda^2$, so

$$g_l(2) \ > \ \frac{(1+4\lambda+4\lambda^2)(1-2\lambda^2)}{(1+4\lambda-12\lambda^2)(1-50\lambda^2+769\lambda^4)} \ = \ \left(1+\frac{16\lambda^2}{1+4\lambda-12\lambda^2}\right)\left(1+\frac{48\lambda^2-769\lambda^4}{1-50\lambda^2+769\lambda^4}\right)$$

$$> \left(1+\frac{16\lambda^2}{1+4\lambda}\right)\left[1+\frac{48\lambda^2\{1-(769/48)\lambda^2\}}{1-50\lambda^2+769\lambda^4}\right] > \left(1+\frac{16\lambda^2}{1+4\lambda}\right)\left(1+48\lambda^2\right)$$

since for $\lambda \le 1/6$, $\dfrac{1-(769/48)\lambda^2}{1-50\lambda^2+769\lambda^4} > 1$ if and only if $50-769/48 > 769\lambda^2$ or $\lambda < \left(\dfrac{50}{769}-\dfrac{1}{48}\right) \approx 0.21$.

As well, $\dfrac{1}{1+4\lambda} > 1-4\lambda$, so

$$g_l(2) \ > \ \{1+16\lambda^2(1-4\lambda)\}(1+48\lambda^2) \ = \ 1+64\lambda^2-64\lambda^3+768\lambda^4-3072\lambda^5$$

$$= \ 1+64\lambda^2(1-\lambda+12\lambda^2-48\lambda^3),$$

and since $12\lambda^2 > 48\lambda^3$ for $\lambda < 1/6$, we have $g_l(2) \ > \ 1+64\lambda^2(1-\lambda)$. But for $\lambda \le 1/6$, $64(1-\lambda) \ge (5/6)64 = 160/3$, so $g_l(2) \ > \ 1+(160/3)\lambda^2$, which is (b).

This concludes the proof of assertion 1.1(v) and Theorem 1.1. □

Next we turn to showing that the largest value of *j* such that $d_l(j) < 0$ is given by $a_l = \lfloor \sqrt{l}/2 \rfloor$.

**3. Proof of Theorem 1.2.**

Proof of assertion 1.2(i). $g_{4a^2}(a) = \dfrac{(4a^2+a)^{(2a)}(4a^2+a)^{(2a)}}{(4a^2+3a)^{(2a)}(4a^2-a)^{(2a)}} = \prod_{i=0}^{2a-1} \dfrac{(4a^2+a-i)^2}{(4a^2+3a-i)(4a^2-a-i)}$

$$= \prod_{i=0}^{2a-1} \frac{\{1+(\alpha/4)(1-i\alpha)\}^2}{\{1+(\alpha/4)(3-i\alpha)\}\{1-(\alpha/4)(1+i\alpha)\}}, \text{ where } \alpha = a^{-1}, \text{ so}$$

$$\log g_{4a^2}(a) = \sum_{i=0}^{2a-1} \left[ 2\log\{1+(\alpha/4)(1-i\alpha)\} - \log\{1+(\alpha/4)(3-i\alpha)\} - \log\{1-(\alpha/4)(1+i\alpha)\} \right]. \qquad (3.1)$$

We approximate the three terms in each summand of (3.1) with the first five powers of $u$ in the log expansions $\log(1+u) = u - u^2/2 + u^3/3 - + \cdots$ and $-\log(1-u) = u + u^2/2 + u^3/3 + \cdots$, say $T_{1i}, \ldots, T_{5i}$, and then summing over $i$, we'll retain the first four powers of $\alpha$ for the expansion of $\log g_{4a^2}(a)$.

First, $T_{1i} = (\alpha/4)\{2(1-i\alpha) - (3-i\alpha) + (1+i\alpha)\} = 0$.

Second, $T_{2i} = (1/2)(\alpha/4)^2\{-2(1-i\alpha)^2 + (3-i\alpha)^2 + (1+i\alpha)^2\}$, where the term in braces equals $-2+4i\alpha-2i^2\alpha^2+9-6i\alpha+i^2\alpha^2+1+2i\alpha+i^2\alpha^2 \;=\; 8$, so $T_{2i} = \alpha^2/4$. Then summing over $i$, we obtain one first-order term in $\alpha$ and no others, which we'll write as $T_2$:

$$T_2 = \sum_{i=0}^{2a-1} \alpha^2/4 = (2a)(\alpha^2/4) = \alpha/2.$$

Third, $T_{3i} = (1/3)(\alpha/4)^3\{2(1-i\alpha)^3 - (3-i\alpha)^3 + (1+i\alpha)^3\}$, where the term in braces equals $(2-6i\alpha+6i^2\alpha^2-2i^3\alpha^3) - (27-3\cdot 9i\alpha+3\cdot 3i^2\alpha^2-i^3\alpha^3) + (1+3i\alpha+3i^2\alpha^2+i^3\alpha^3) = -24(1-i\alpha)$, so $T_{3i} = -(\alpha^3/8)(1-i\alpha)$. Then summing over $i$, we obtain a third-order term in $\alpha$ and no others,

$$T_3 = \sum_{i=0}^{2a-1} T_{3i} = -(2a)(\alpha^3/8) + \frac{2a(2a-1)}{2}\alpha^4/8 = -\alpha^2/4 + \alpha^2/4 - \alpha^3/8 = -\alpha^3/8.$$

Fourth, $T_{4i} = (1/4)(\alpha/4)^4\{-2(1-i\alpha)^4 + (3-i\alpha)^4 + (1+i\alpha)^4\}$, where the term in braces equals

$$(-2+8i\alpha-12i^2\alpha^2+8i^3\alpha^3-2i^4\alpha^4) + (81-4\cdot 27i\alpha+6\cdot 9i^2\alpha^2-4\cdot 3i^3\alpha^3+i^4\alpha^4)$$
$$+ (1+4i\alpha+6i^2\alpha^2+4i^3\alpha^3+i^4\alpha^4) \;=\; 80-96i\alpha+48i^2\alpha^2 = 16(5-6i\alpha+3i^2\alpha^2),$$

so $T_{4i} = (\alpha^4/64)(5-6i\alpha+3i^2\alpha^2)$. Then summing over $i$, we get a sole other third-order term in $\alpha$,

$$T_4 = \sum_{i=0}^{2a-1} T_{4i} = (2a)(\frac{5}{64})\alpha^4 + \frac{2a(2a-1)}{2}(\frac{-6}{64})\alpha^5 + \frac{2a(2a-1)(4a-1)}{6}(\frac{3}{64})\alpha^6$$

$$= \frac{5}{32}\alpha^3 - \frac{3}{16}\alpha^3 + \frac{3}{32}\alpha^4 + \frac{1}{8}\alpha^3 - \frac{3}{32}\alpha^4 + \alpha^5/64 = (\frac{5-6+4}{32})\alpha^3 + \alpha^5/64$$

$$= \frac{3}{32}\alpha^3 + O(\alpha^5)\,.$$

Fifth, $T_{5i} = (1/5)(\alpha/4)^5\{2(1-i\alpha)^5 - (3-i\alpha)^5 + (1+i\alpha)^5\}$, where the term in braces equals

$$(2-10i\alpha+20i^2\alpha^2-20i^3\alpha^3+10i^4\alpha^4-2i^5\alpha^5) - (243-5\cdot 81i\alpha+10\cdot 27i^2\alpha^2-10\cdot 9i^3\alpha^3+5\cdot 3i^4\alpha^4-i^5\alpha^5)$$

$$+(1+5i\alpha+10i^2\alpha^2+10i^3\alpha^3+5i^4\alpha^4+i^5\alpha^5) \;=\; -240+400i\alpha-240i^2\alpha^2+80i^3\alpha^3$$

$$= 80(-3+5i\alpha-3i^2\alpha^2+i^3\alpha^3)\,,$$

so $T_{5i} = (\alpha^5/64)(-3+5i\alpha-3i^2\alpha^2+i^3\alpha^3)$. Then summing over $i$, $T_5 = \sum_{i=0}^{2a-1} T_{5i} =$

$$(2a)(\frac{-3}{64})\alpha^5 + \frac{2a(2a-1)}{2}(\frac{5}{64})\alpha^6 + \frac{2a(2a-1)(4a-1)}{6}(\frac{-3}{64})\alpha^7 + \frac{(2a)^2(2a-1)^2}{4}(\frac{1}{64})\alpha^3$$

and retaining only terms of order $\alpha^4$, its coefficient is $-\frac{3}{32} + \frac{5}{32} - \frac{1}{8} + \frac{1}{16} = \frac{-3+5-4+2}{32} = 0$.

Higher order terms of the form $T_{pi} = p^{-1}(\alpha/4)^p\{c_0 + c_1 i\alpha + \cdots + c_p(i\alpha)^p\}$ only yield terms of order $\alpha^{p-1}$ or higher after summing over $i$, so we may stop here. This establishes (1.22), namely,

$$\log g_{4a^2}(a) = \frac{\alpha}{2} - \frac{\alpha^3}{8} + \frac{3\alpha^3}{32} + O(\alpha^5) = \frac{\alpha}{2} - \frac{\alpha^3}{32} + O(\alpha^5)\,. \tag{3.2}$$

For $h_{4a^2}(a) = \frac{4a^2+a}{4a^2-a} = \frac{1+\alpha/4}{1-\alpha/4}$, we expand

$$\log h_{4a^2}(a) = \alpha/4 - (\alpha/4)^2/2 + (\alpha/4)^3/3 - (\alpha/4)^4/4 + - \cdots$$

$$+\alpha/4 + (\alpha/4)^2/2 + (\alpha/4)^3/3 + (\alpha/4)^4/4 + \cdots$$

$$= \alpha/2 + \alpha^3/96 + O(\alpha^5) > \alpha/2 + \alpha^3/96\,,$$

since all non-zero coefficients are positive. This establishes (1.23) and assertion 1.2(i), which is that ignoring terms of order $a^{-5}$, $g_{4a^2}(a) < h_{4a^2}(a)$.

Proof of assertion 1.2(ii). We can use (2.7) for the ratio $g_l(a+1)/g_l(a)$ to write

$$g_{4a^2}(a+1) = g_{4a^2}(a)\frac{(4a^2+a+1)^3(4a^2-a)^3}{(4a^2+3a+3)^{(3)}(4a^2-3a)^{(3)}},$$

where the second factor equals

$$\frac{\{1+(\alpha/4)(1+\alpha)\}^3\{1-(\alpha/4)\}^3}{\{1+(\alpha/4)(3+3\alpha)\}\,\{1+(\alpha/4)(3+2\alpha)\}\,\{1+(\alpha/4)(3+\alpha)\}\,\{1-3\alpha/4\}\,\{1-(\alpha/4)(3+\alpha)\}\,\{1-(\alpha/4)(3+2\alpha)\}}$$

with logarithm approximated as follows.

$$\begin{array}{lllll}
3(\frac{\alpha}{4})(1+\alpha) & -(3/2)(\frac{\alpha}{4})^2(1+\alpha)^2 & +(3/3)(\frac{\alpha}{4})^3(1+\alpha)^3 & -(3/4)(\frac{\alpha}{4})^4(1+\alpha)^4 & +-\cdots \\
-3(\frac{\alpha}{4}) & -(3/2)(\frac{\alpha}{4})^2 & -(3/3)(\frac{\alpha}{4})^3 & -(3/4)(\frac{\alpha}{4})^4 & -\cdots \\
-3(\frac{\alpha}{4})(1+\alpha) & +(9/2)(\frac{\alpha}{4})^2(1+\alpha)^2 & -(27/3)(\frac{\alpha}{4})^3(1+\alpha)^3 & +(81/4)(\frac{\alpha}{4})^4(1+\alpha)^4 & -+\cdots \\
-(\frac{\alpha}{4})(3+2\alpha) & +(1/2)(\frac{\alpha}{4})^2(3+2\alpha)^2 & -(1/3)(\frac{\alpha}{4})^3(3+2\alpha)^3 & +(1/4)(\frac{\alpha}{4})^4(3+2\alpha)^4 & -+\cdots \\
-(\frac{\alpha}{4})(3+\alpha) & +(1/2)(\frac{\alpha}{4})^2(3+\alpha)^2 & -(1/3)(\frac{\alpha}{4})^3(3+\alpha)^3 & +(1/4)(\frac{\alpha}{4})^4(3+\alpha)^4 & -+\cdots \\
+3(\frac{\alpha}{4}) & +(9/2)(\frac{\alpha}{4})^2 & +(27/3)(\frac{\alpha}{4})^3 & +(81/4)(\frac{\alpha}{4})^4 & -\cdots \\
+(\frac{\alpha}{4})(3+\alpha) & +(1/2)(\frac{\alpha}{4})^2(3+\alpha)^2 & +(1/3)(\frac{\alpha}{4})^3(3+\alpha)^3 & +(1/4)(\frac{\alpha}{4})^4(3+\alpha)^4 & +\cdots \\
+(\alpha/4)(3+2\alpha) & +(1/2)(\alpha/4)^2(3+2\alpha)^2 & +(1/3)(\alpha/4)^3(3+2\alpha)^3 & +(1/4)(\alpha/4)^4(3+2\alpha)^4 & +\cdots.
\end{array}$$

Collecting the coefficients of the first four powers of $\alpha$, we find

$\text{Coeff}(\alpha) = (1/4)(3-3-3-3-3+3+3+3) = 0$.

$\text{Coeff}(\alpha^2) = (1/4)(3+0-3-2-1+0+1+2)$ from the 1st expansion term

$\qquad +(1/32)(-3-3+9+9+9+9+9+9)$ " " 2nd " "

$\qquad = 48/32 = 3/2$.

$\text{Coeff}(\alpha^3) = (1/16)(-3+0+9+6+3+0+3+6)$ from the 2nd expansion term

$\qquad +(1/64)(1-1-9-9-9+9+9+9)$ " " 3rd " "

$\qquad = 24/16 = 3/2$.

$\text{Coeff}(\alpha^4) = (1/32)(-3+0+9+4+1+0+1+4)$ from the 2nd expansion term

$\qquad +(1/64)(3+0-27-18-9+0+9+18)$ " " 3rd " "

$\qquad +(1/1024)(-3-3+81+81+81+81+81+81)$ " " 4th " "

$$= \frac{16}{32} - \frac{24}{64} + \frac{480}{1024} = \frac{16-12+15}{32} = \frac{19}{32}.$$

Thus, the log of the second factor is $\frac{3}{2}\alpha^2 + \frac{3}{2}\alpha^3 + \frac{19}{32}\alpha^4 + O(\alpha^5)$, and combining this with (3.2) yields

$$\log g_{4a^2}(a+1) = \left(\frac{\alpha}{2} - \frac{\alpha^3}{32}\right) + \left(\frac{3}{2}\alpha^2 + \frac{3}{2}\alpha^3 + \frac{19}{32}\alpha^4\right) + O(\alpha^5) = \frac{\alpha}{2} + \frac{3}{2}\alpha^2 + \frac{47\alpha^3}{32} + \frac{19\alpha^4}{32} + O(\alpha^5), \quad (3.3)$$

which is (1.24).

For $h_{4a^2}(a+1) = \dfrac{4a^2+a+1}{4a^2-a-1} = \dfrac{1+(\alpha/4)(1+\alpha)}{1-(\alpha/4)(1+\alpha)}$,

$$\log h_{4a^2}(a+1) = (\alpha/4)(1+\alpha) - (\alpha/4)^2(1+\alpha)^2/2 + (\alpha/4)^3(1+\alpha)^3/3 - (\alpha/4)^4(1+\alpha)^4/4 + - \cdots$$

$$+(\alpha/4)(1+\alpha) + (\alpha/4)^2(1+\alpha)^2/2 + (\alpha/4)^3(1+\alpha)^3/3 + (\alpha/4)^4(1+\alpha)^4/4 + \cdots$$

$$= (\frac{\alpha}{2})(1+\alpha) + \frac{2}{3}(\alpha/4)^3(1+\alpha)^3 + O(\alpha^5) = \frac{\alpha}{2} + \frac{\alpha^2}{2} + \frac{\alpha^3}{96} + \frac{\alpha^4}{32} + O(\alpha^5),$$

and comparing this with (3.3) establishes (1.25) and assertion 1.2(ii), which is that ignoring terms of order $a^{-5}$, $g_{4a^2}(a+1) > h_{4a^2}(a+1)$.

<u>Proof of assertion 1.2(iii)</u>.

$$g_{4a^2-1}(a) = \frac{(4a^2+a-1)^{(2a)}(4a^2+a-1)^{(2a)}}{(4a^2+3a-1)^{(2a)}(4a^2-a-1)^{(2a)}} = \prod_{i=0}^{2a-1} \frac{(4a^2+a-1-i)^2}{(4a^2+3a-1-i)(4a^2-a-1-i)}$$

$$= \prod_{i=0}^{2a-1} \frac{\left[1+(\alpha/4)\{1-\alpha(i+1)\}\right]^2}{\left[1+(\alpha/4)\{3-\alpha(i+1)\}\right]\left[1-(\alpha/4)\{1+\alpha(i+1)\}\right]} = \prod_{i=1}^{2a} \frac{\{1+(\alpha/4)(1-i\alpha)\}^2}{\{1+(\alpha/4)(3-i\alpha)\}\{1-(\alpha/4)(1+i\alpha)\}}, \text{ so}$$

$$\log g_{4a^2-1}(a) = \sum_{i=1}^{2a} \left[2\log\{1+(\alpha/4)(1-i\alpha)\} - \log\{1+(\alpha/4)(3-i\alpha)\} - \log\{1-(\alpha/4)(1+i\alpha)\}\right]$$

which is the same as (3.1) except for the limits of summation. Thus we have $T_1 = 0$ and $T_2 = \alpha/2$ as before. Summing $T_{3i} = -(\alpha^3/8)(1-i\alpha)$ from above over $i = 1,\ldots,2a$, we have

$$T_3 = \sum_{i=1}^{2a} T_{3i} = -(2a)(\alpha^3/8) + \frac{2a(2a+1)}{2}\alpha^4/8 = -(\alpha^2/4) + (\alpha^2/4) + (\alpha^3/8) = \alpha^3/8,$$

and summing $T_{4i} = (\alpha^4/64)(5-6i\alpha+3i^2\alpha^2)$ from above over $i=1,\ldots,2a$, we have

$$T_4 = \sum_{i=1}^{2a} T_{4i} = (2a)(\frac{5}{64})\alpha^4 + \frac{2a(2a+1)}{2}(\frac{-6}{64})\alpha^5 + \frac{2a(2a+1)(4a+1)}{6}(\frac{3}{64})\alpha^6$$

$$= \frac{5}{32}\alpha^3 - \frac{3}{16}\alpha^3 - \frac{3}{32}\alpha^4 + \frac{1}{8}\alpha^3 + \frac{3}{32}\alpha^4 + +O(\alpha^5) = (\frac{5-6+4}{32})\alpha^3 + O(\alpha^5) = \frac{3}{32}\alpha^3 + O(\alpha^5).$$

Finally, $T_{5i} = (\alpha^5/64)(-3+5i\alpha-3i^2\alpha^2+i^3\alpha^3)$ from above, so summing over $i=1,\ldots,2a$,

$$T_5 = (2a)(\frac{-3}{64})\alpha^5 + \frac{2a(2a+1)}{2}(\frac{5}{64})\alpha^6 + \frac{2a(2a+1)(4a+1)}{6}(\frac{-3}{64})\alpha^7 + \frac{(2a)^2(2a+1)^2}{4}(\frac{1}{64})\alpha^8.$$

Retaining only terms of order $\alpha^4$, its coefficient is $-\frac{3}{32}+\frac{5}{32}-\frac{1}{8}+\frac{1}{16} = \frac{-3+5-4+2}{32} = 0$ as before.

Thus,

$$\log g_{4a^2-1}(a) = \frac{\alpha}{2} + \frac{\alpha^3}{8} + \frac{3\alpha^3}{32} + O(\alpha^5) = \frac{\alpha}{2} + \frac{7\alpha^3}{32} + O(\alpha^5) \tag{3.4}$$

which establishes (1.26).

For $h_{4a^2-1}(a)$, we have $h_{4a^2-1}(a) = \frac{4a^2+a-1}{4a^2-a-1} = \frac{1+(\alpha/4)(1-\alpha)}{1-(\alpha/4)(1+\alpha)}$, so that

$$\log h_{4a^2-1}(a) = (\alpha/4)(1-\alpha) - \frac{1}{2}(\alpha/4)^2(1-\alpha)^2 + \frac{1}{3}(\alpha/4)^3(1-\alpha)^3 - \frac{1}{4}(\alpha/4)^4(1-\alpha)^4 + - \cdots$$

$$+ (\alpha/4)(1+\alpha) + \frac{1}{2}(\alpha/4)^2(1+\alpha)^2 + \frac{1}{3}(\alpha/4)^3(1+\alpha)^3 + \frac{1}{4}(\alpha/4)^4(1+\alpha)^4 + \cdots$$

$$= \frac{\alpha}{2} + \frac{\alpha^3}{8} + \frac{\alpha^3}{96} + O(\alpha^5) = \frac{\alpha}{2} + \frac{13\alpha^3}{96} + O(\alpha^5) > \frac{\alpha}{2} + \frac{13\alpha^3}{96},$$

which establishes (1.27) and assertion 1.2(iii), that ignoring terms of order $a^{-5}$, $g_{4a^2-1}(a) > h_{4a^2-1}(a)$.

Proof of assertion 1.2(iv). We prove something slightly more general, namely that for given $a \ge 3$ and each $j$ with $a \le j \le 2a+1$, $g_l(j)/h_l(j)$ decreases strictly in $l$ over the interval $4a^2 \le l \le 4(a+1)^2-1$. See Remark 3.1 below for why we need the restriction $a \ge 3$. Let $l = x + 4a^2$ for continuous $x$ satisfying $0 \le x \le 4(a+1)^2 - 1 - 4a^2 = 8a+3$. Define the log ratio $L(x) = L(x, j)$ by

$$\exp L(x) := \frac{g_{x+4a^2}(j)}{h_{x+4a^2}(j)} = \frac{(x+4a^2+j)^{(2j)}(x+4a^2+j)^{(2j)}}{(x+4a^2+3j)^{(2j)}(x+4a^2-j)^{(2j)}} \cdot \frac{x+4a^2-j}{x+4a^2+j}$$

$$= \prod_{i=0}^{2j-1} \left\{ \frac{(x+4a^2+j-i)^2}{(x+4a^2+3j-i)(x+4a^2-j-i)} \right\} \cdot \frac{x+4a^2-j}{x+4a^2+j}$$

so that the derivative of $L(x)$ with respect to $x$ is

$$L'(x)=\sum_{i=0}^{2j-1}\left\{\frac{2}{x+4a^2+j-i}-\frac{1}{x+4a^2+3j-i}-\frac{1}{x+4a^2-j-i}\right\}+\frac{1}{x+4a^2-j}-\frac{1}{x+4a^2+j}$$
$$=\sum_{i=0}^{2j-1}\left\{\frac{2}{x_i+j}-\frac{1}{x_i+3j}-\frac{1}{x_i-j}\right\}+\frac{1}{x_0-j}-\frac{1}{x_0+j}, \tag{3.5}$$

where $x_i=x+4a^2-i$. Now $\frac{1}{x_0-j}-\frac{1}{x_0+j}=\frac{2j}{x_0^2-j^2}$ and the summands in (3.5) are the negatives of

$$\frac{1}{x_i+3j}+\frac{1}{x_i-j}-\frac{2}{x_i+j}=\frac{1}{x_i+3j}-\frac{x_i-3j}{x_i^2-j^2}=\frac{(x_i^2-j^2)-(x_i^2-9j^2)}{(x_i+3j)(x_i^2-j^2)}=\frac{8j^2}{(x_i+3j)(x_i^2-j^2)}.$$

Thus $L'(x)<0$ if and only if $\sum_{i=0}^{2j-1}\frac{8j^2}{(x_i+3j)(x_i^2-j^2)}>\frac{2j}{x_0^2-j^2}$ or

$$\sum_{i=0}^{2j-1}\frac{4j}{x_i+3j}\left\{\frac{x_0^2-j^2}{x_i^2-j^2}\right\}>1. \tag{3.6}$$

But $x_0>x_1>\cdots$, so the left-hand side of (3.6) exceeds $(2j)\frac{4j}{x_0+3j}=\frac{8j^2}{x+4a^2+3j}$ and since $x\le 8a+3$, the left-hand side of (3.6) exceeds $\frac{8j^2}{(8a+3)+4a^2+3j}$. Thus it suffices to show that $8j^2\ge 4a^2+8a+3(j+1)$. Now for $j=a$, this condition is $4a^2-11a-3\ge 0$, which holds for $a\ge 3$. When $a+1\le j\le 2a+1$, $8j^2\ge 8(a+1)^2=8a^2+16a+8$ and $4a^2+8a+3(j+1)\le 4a^2+14a+6$, so it suffices to show that $8a^2+16a+8\ge 4a^2+14a+6$, or $4a^2+2a+2\ge 0$, which is true for all $a$. Thus $L(x)$ decreases in $x$, whence $g_l(j)/h_l(j)$ decreases in $l$ for $4a^2\le l\le 4(a+1)^2-1$, as asserted. Assertion 1.2(iv) is the case of $j=a$ and $j=a+1$.

This concludes the proof of Theorem 1.2. □

Remark 3.1. The restriction to $a\ge 3$ in Theorem 1.2(iv) is needed because for $a=1=j$, $g_l(1)/h_l(1)$ actually increases over $l=7,\ldots,15$ and for $a=2=j$, $g_l(2)/h_l(2)$ increases over $l=31,\ldots,35$. Nevertheless, as indicated in Section 1, the formula $a_l=\lfloor\sqrt{l}/2\rfloor$ holds for $a=1$ and 2 as well, since direct calculation shows that $g_l(1)<h_l(1)$ while $g_l(2)>h_l(2)$ for $l=6,\ldots,15$, and $g_l(2)<h_l(2)$ while $g_l(3)>h_l(3)$ for $l=16,\ldots,35$. □

## 4. Proof of Theorem 1.3.

Given any integer $a \geq 1$, let $l \geq 9$ be such that $4a^2 \leq l \leq 4(a+1)^2 - 1$, so that $a_l = a$ for such *l*. For convenience in this section we will suppress the subscript *l* on $g_l(j)$, $h_l(j)$, $d_l(j)$, and $Q_l(j)$ when no confusion would arise. It will also be very helpful to re-express the pairwise inequalities for *d*(*j*)*Q*(*j*) in (1.29) as follows. For each $a$ and $j \in \{1,\dots,a\}$, multiply both $d(j+a+1)/\{-d(j)\}$ and $Q(j)/Q(j+a+1)$ by the ratio $h(j)/h(j+a+1)$ to obtain

$$\frac{h(j)}{h(j+a+1)}\frac{d(j+a+1)}{-d(j)} = \frac{h(j)}{h(j+a+1)}\frac{g(j+a+1)-h(j+a+1)}{h(j)-g(j)} = \frac{\dfrac{g(j+a+1)}{h(j+a+1)}-1}{1-\dfrac{g(j)}{h(j)}} \quad \text{and}$$

$$\frac{h(j)}{h(j+a+1)}\frac{Q(j)}{Q(j+a+1)} = \frac{j(j+l+a+1)}{(j+l)(j+a+1)}\frac{\dbinom{2l-1+2j}{l-1+j}\dbinom{2l-1-2j}{l-1+j}}{\dbinom{2l+1+2a+2j}{l+a+j}\dbinom{2l-3-2a-2j}{l+a+j}},$$

since from (1.28), the product of $h(j)/h(j+a+1)$ with the leading factors in $Q(j)/Q(j+a+1)$ is

$$\frac{h(j)}{h(j+a+1)}\cdot\frac{j(l-j)}{(l+j)^2}\frac{(l+j+a+1)^2}{(j+a+1)(l-j-a-1)} = \frac{(l+j)}{(l-j)}\frac{(l-j-a-1)}{(l+j+a+1)}\cdot\frac{j(l-j)}{(l+j)^2}\frac{(l+j+a+1)^2}{(j+a+1)(l-j-a-1)}$$
$$= \frac{j(j+l+a+1)}{(j+l)(j+a+1)}.$$

Then (1.29), or equivalently $d(j+a+1)/\{-d(j)\} > Q(j)/Q(j+a+1)$, holds if and only if

$$\frac{g(j+a+1)}{h(j+a+1)} - 1 > \left\{1-\frac{g(j)}{h(j)}\right\}\left\{\frac{j(j+l+a+1)}{(j+l)(j+a+1)}\right\}\frac{\dbinom{2l-1+2j}{l-1+j}\dbinom{2l-1-2j}{l-1+j}}{\dbinom{2l+1+2a+2j}{l+a+j}\dbinom{2l-3-2a-2j}{l+a+j}}. \tag{4.1}$$

Some remarkable monotonicities appear with re-expression (4.1). For example, Figure 4.1 plots values of $d(j+a+1)Q(j+a+1)$ (red boxes) lying above the values of $-d(j)Q(j)$ (blue diamonds) against $j$=1,…,8 in the case $l$=256 with $a_l = 8$, illustrating the pairwise inequalities. Note how the vertical differences between corresponding points are not monotonic in *j*. Figure 4.2 shows the values on the left-hand side of (4.1) (red boxes) also lying above the values on the right-hand side of (4.1) (blue diamonds) in the same case. Note how the vertical differences now strictly increase in *j*. Figures 4.3 to 4.5 reveal another monotonicity, this time in *l*. Figure 4.3 plots the left- and right-hand sides of (4.1) against $l$=36,…,63 in the case $a$=3, $j$=1. Figure 4.4 does the same for $j$=2 and Figure 4.5 does the same for $j$=3. Note how in each case the vertical differences between corresponding points strictly *decrease* in *l*, even though the blue diamonds in the latter case are not themselves monotone.

These monotonicities appear to hold in general, which suggests that to prove (4.1) it would suffice to prove it for each $a$ at the largest $l=4(a+1)^2-1$ and the smallest $j=1$. We call this the *most stringent case*. Yet another monotonicity appears in Table 4.1 below. Let $\Delta(a,j,l)$ denote the differences between left- and right-hand sides of (4.1) (not the forward difference operator here). Table 4.1 shows values of $\Delta(a,j,4(a+1)^2-1)$ for $a=1,\ldots,10$. Note how these minimum values of $\Delta(a,j,l)$ over $l$ decrease in $a$ for each $j$, and in particular for the values $\Delta(a,1,4(a+1)^2-1)$ in the most stringent case. See also Remark 4.1 below.

Table 4.1
Differences between the left- and right-hand sides of (4.1) evaluated at $l=4(a+1)^2-1$.

| $a$ | $j=1$ | $j=2$ | $j=3$ | $j=4$ | $j=5$ | $j=6$ | $j=7$ | $j=8$ | $j=9$ | $j=10$ |
|---|---|---|---|---|---|---|---|---|---|---|
| 1 | 0.440909 | | | | | | | | | |
| 2 | 0.134941 | 0.360421 | | | | | | | | |
| 3 | 0.069596 | 0.144741 | 0.325816 | | | | | | | |
| 4 | 0.043015 | 0.085709 | 0.142530 | 0.300135 | | | | | | |
| 5 | 0.029335 | 0.058016 | 0.088660 | 0.140139 | 0.277542 | | | | | |
| 6 | 0.021317 | 0.042164 | 0.062874 | 0.088736 | 0.137890 | 0.257280 | | | | |
| 7 | 0.016203 | 0.032102 | 0.047553 | 0.064230 | 0.088335 | 0.135406 | 0.239180 | | | |
| 8 | 0.012736 | 0.025280 | 0.037415 | 0.049678 | 0.064481 | 0.087874 | 0.132574 | 0.223069 | | |
| 9 | 0.010276 | 0.020430 | 0.030268 | 0.039933 | 0.050440 | 0.064425 | 0.087350 | 0.129440 | 0.208735 | |
| 10 | 0.008467 | 0.016856 | 0.025011 | 0.032935 | 0.041072 | 0.050665 | 0.064294 | 0.086710 | 0.126098 | 0.195960 |

We conjecture these properties hold in general, as follows.

Conjecture 4.1 [monotonicity properties (MP)]

(i) For each $a\geq 1$ and $l\geq 9$ in the discrete interval $I(a):=\{4a^2,\ldots,4(a+1)^2-1\}$, $\Delta(a,j,l)$ increases strictly in $j=1,\ldots,a$, so that $\arg\min_j \Delta(a,j,l) = 1$ for each such $a$ and $l$.

(ii) For each $a\geq 1$ and $j=1,\ldots,a$, $\Delta(a,j,l)$ decreases strictly in $l$ over $I(a)\cap\{l\geq 9\}$, so that $\arg\min_l \Delta(a,j,l) = 4(a+1)^2-1$ for each such $a$ and $j$.

(iii) For each $j\geq 1$, $\Delta(a,j,4(a+1)^2-1)=\min_l \Delta(a,j,l)$ decreases strictly in $a$ for $a\geq j$.

In particular, $\min_{j,l}\Delta(a,j,l)=\Delta(a,1,4(a+1)^2-1)$ decreases strictly in $a$ for $a\geq 1$. □

We will assume MP going forward and proceed to show that $\Delta(a,1,4(a+1)^2-1)>0$ for sufficiently large $a$. That will establish the pairwise inequalities for all $a$, $j$, and $l$ under MP.

Let *LHS* denote the left-hand side of (4.1) and *RHS* the right-hand side, each evaluated at $j=1$ and $l=4(a+1)^2-1$ in the most stringent case. *Going forward it will be convenient to shift a by one* *unit so that we may write* $4a^2-1$ *instead of* $4(a+1)^2-1$. To that end, let $a'=a+1$, in which case the left-hand side of (4.1) becomes $LHS=\{g_{4a'^2-1}(a'+1)/h_{4a'^2-1}(a'+1)\}-1$. Then, dropping the primes and remembering that the original $a$ is now $a-1$ and the original $j+a+1$ is now $a+1$, we have

$$
\begin{aligned}
LHS &= \frac{(4a^2-1+a+1)^{(2a+2)}(4a^2-1+a+1)^{(2a+2)}(4a^2-1-a-1)}{(4a^2-1+3a+3)^{(2a+2)}(4a^2-1-a-1)^{(2a+2)}(4a^2-1+a+1)}-1 \\
&= \frac{(4a^2+a-1)^{(2a+1)}(4a^2+a)^{(2a+2)}}{(4a^2-a-3)^{(2a+1)}(4a^2+3a+2)^{(2a+2)}}-1 .
\end{aligned}
\tag{4.2}
$$

Similarly, after shifting $a$ by 1 and dropping primes,

$$
RHS=\left\{1-\frac{g_{4a^2-1}(1)}{h_{4a^2-1}(1)}\right\}\left\{\frac{1+(4a^2-1)+a)}{(1+4a^2-1)(a+1)}\right\}\frac{\binom{2(4a^2-1)-1+2}{(4a^2-1)-1+1}\binom{2(4a^2-1)-1-2}{(4a^2-1)-1+1}}{\binom{2(4a^2-1)+1+2(a-1)+2}{(4a^2-1)+(a-1)+1}\binom{2(4a^2-1)-3-2(a-1)-2}{(4a^2-1)+(a-1)+1}}
$$

$$
=\left\{1-\frac{(4a^2)^{(2)}(4a^2)^{(2)}(4a^2-1-1)}{(4a^2+2)^{(2)}(4a^2-2)^{(2)}(4a^2-1+1)}\right\}\left\{\frac{4a+1}{4a(a+1)}\right\}\frac{\binom{8a^2-1}{4a^2-1}\binom{8a^2-5}{4a^2-1}}{\binom{8a^2+2a-1}{4a^2+a-1}\binom{8a^2-2a-5}{4a^2+a-1}}
$$

$$
=\left\{1-\frac{(4a^2)^2(4a^2-1)^2(4a^2-2)}{(4a^2+2)(4a^2+1)(4a^2-2)(4a^2-3)(4a^2)}\right\}\left\{\frac{4a+1}{4a(a+1)}\right\}\frac{\binom{8a^2-1}{4a^2-1}\binom{8a^2-5}{4a^2-1}}{\binom{8a^2+2a-1}{4a^2+a-1}\binom{8a^2-2a-5}{4a^2+a-1}}
$$

$$
=\left\{1-\frac{4a^2(4a^2-1)^2}{(4a^2+2)(4a^2+1)(4a^2-3)}\right\}\left\{\frac{4a+1}{4a(a+1)}\right\}\frac{\binom{8a^2-1}{4a^2-1}\binom{8a^2-5}{4a^2-1}}{\binom{8a^2+2a-1}{4a^2+a-1}\binom{8a^2-2a-5}{4a^2+a-1}} .
$$

But the first factor in braces equals

$$
\frac{(16a^4+12a^2+2)(4a^2-3)-4a^2(16a^4-8a^2+1)}{(4a^2+2)(4a^2+1)(4a^2-3)}=\frac{32a^2(a^2-1)-6}{(4a^2+2)(4a^2+1)(4a^2-3)} ,
$$

so

$$
RHS=\left[\frac{\{32a^2(a^2-1)-6\}(4a+1)/(4a)}{(4a^2+2)(4a^2+1)(4a^2-3)(a+1)}\right]\frac{\binom{8a^2-1}{4a^2-1}\binom{8a^2-5}{4a^2-1}}{\binom{8a^2+2a-1}{4a^2+a-1}\binom{8a^2-2a-5}{4a^2+a-1}} .
\tag{4.3}
$$

Lemma 4.1. *Let* $\alpha = a^{-1}$. *Then as* $a \to \infty$ *and* $\alpha \to 0$, $\log(1+LHS) = \alpha^2 + (19/12)\alpha^3 + o(\alpha^3)$. *Thus* $LHS = \alpha^2 + (19/12)\alpha^3 + o(\alpha^3)$ *as well.*

Proof. See Appendix C. □

Lemma 4.2. *As* $a \to \infty$ *and* $\alpha \to 0$, $RHS = (e/2)\alpha^3 + o(\alpha^3)$.

Proof. See Appendix D. □

It follows from Lemmas 4.1 and 4.2 that the difference between left- and right-hand sides of (4.1), now written with shifted $a$ as $\Delta(a-1, 1, 4a^2-1)$, equals $\alpha^2 + (\frac{19}{12} - \frac{e}{2})\alpha^3 + o(\alpha^3)$, or dividing by $\alpha^2$, that

$$\Delta(a-1, 1, 4a^2-1)/\alpha^2 = 1 + (\frac{19}{12} - \frac{e}{2})\alpha\{1 + o(1)\},$$

which exceeds 1 for sufficiently small $\alpha$, since 19/12 > $e$/2. Thus $\Delta(a-1, 1, 4a^2-1) > \alpha^2 > 0$ for sufficiently large $a$, *which establishes the pairwise inequalities for* $d(j)Q(j)$ *of Theorem 1.3 under MP.* □

For example, the last line of Table 4.1, which corresponds to (shifted) $a$=11, has $\Delta(10, 1, 483) = 0.008467 > 1/11^2 = 0.008264$, with $\Delta(10, 1, 483)/\alpha^2 = 1.0245$. Figure 4.6 plots $\Delta(a, 1, 4(a+1)^2 - 1)$ and $\alpha^2 + (\frac{19}{12} - \frac{e}{2})\alpha^3 = \frac{1}{(a+1)^2} + (\frac{19}{12} - \frac{e}{2})\left(\frac{1}{a+1}\right)^3$ against (original) $a$=1,…,10, showing good agreement.

Remark 4.1. Figure 4.7 illustrates MP(i) and MP(ii) in the case $a$=5 and $a$=8. These plots display convexity in $l$ for each $j$=1,…,5. *Indeed, numerical evidence suggests that for* $a \geq 4$, $\Delta(a, j, l)$ *is uniformly convex in l for each j=1,…,a.* However, $\Delta(2, j, l)$ is *not* globally convex in $l$ for $j$=1 nor is $\Delta(3, j, l)$ globally convex in $l$ for $j$=2. With respect to convexity in $j$, for $a \leq 4$, $\Delta(a, j, l)$ is globally convex in $j$ for each $l \in I(a)$ (vacuously so for $a$=1 or 2). However, though it is not easy to see in Figure 4.7a, $\Delta(5, j, l)$ is *not* globally convex in $j$ for the first 28 values of $l$=100,…,127, and for $a$=8 in Figure 4.7b (or for *any* $a \geq 6$), $\Delta(a, j, l)$ is not globally convex in $j$ *for any value* of $l \in I(a)$, having instead an ogival shape (concave increasing turning to convex increasing) for each $l$. Of course, MP does not depend on any convexity properties of $\Delta(a, j, l)$. □

## 5. Discussion

There are several gaps to fill for a complete proof that $H_0(n) \ge H_{-1}(n)$ for all $n \ge 1$. The monotonicity properties MP (Conjecture 4.1) is one from the present work and another is the case $r = n \bmod 4 > 0$. In the latter case, the p.f. appearing immediately above (1.10), call it $Q_n(j)$, now puts non-zero mass at $j=0$, and as mentioned in Section 1, in this case $h_n(j) := \{(l+j)/(l-j)\}\{(j-r/4)/(j+r/4)\}$ so that $h_n(0) = -1$ and $d_n(0) = 2$. Now, equation (6.14) of L25 shows that sufficient condition (1.7), which is $E^*_\kappa[K^{-1}] \ge 4/n$ [where the expectation is taken with respect to p.f. (1.6) for $K \in \{1,\ldots,\kappa\}$], holds if $E_n[d_n(J)] > Q_n(0)$ [where the expectation is taken wih respect to $Q_n(j)$ for $J \in \{0,1,\ldots,\tau\}$]. Whereas the pairwise property for $d_l(j)Q_l(j)$ was instrumental for establishing the latter inequality when $n$ is a multiple of 4, when $n \bmod 4 > 0$ it appears that one only needs to consider values of $j$ smaller than the value for which $Q_n(j)$ becomes positive for the last time, *with no consideration needed of the larger pair members*. Specifically, L25 conjectures the following property holds which clearly implies $E_n[d_n(J)] > Q_n(0)$.

Property 6.5 of L25. For $n = 4l + r$ with $l \ge 6$ and $r \in \{1, 2, 3\}$, $d_n(j)$ continues to have strictly increasing first differences, so that $a_{l,r} := \max\{j : d_n(j) \le 0\}$ is well-defined [though not with the simple formula (1.14)]. Then $a_{l,r}$ satisfies $Q_n(0) + E_n[d_n(J) \mid 0 < J \le a_{l,r}] \cdot Q_n[0 < J \le a_{l,r}] \; > \; 0$. □

Thus monotonicity properties MP and Property 6.5 of L25 establish that $H_0(n) > H_{-1}(n)$ for all $n \ge 7$ and $H_0(n) \ge H_{-1}(n)$ for all $n \ge 1$. Elementary arguments then show that $H_0(n) > H_1(n)$ for all $n \ge 7$ and $H_0(n) \ge H_1(n)$ for all $n \ge 1$ (see Corollary 1.1 of L25).

There are other gaps to fill before we have a complete proof of the global unimodality of $H_s(n)$ at $s=0$. To that end, L25 states some additional properties, all of which we conjecture hold true. For establishing that $H_s(n) \ge H_{s+1}(n)$ for each positive $s \le n-1$ [where $H_s(n) = 0$ for $s$ larger than $n-1$], see monotonicity Properties 4.1 and 4.2 or, alternatively, distributional Properties 5.1 to 5.3 of L25. For establishing that $H_s(n) \ge H_{s-1}(n)$ for each negative $s \ge -\lfloor (n-1)/2 \rfloor$ [where $H_s(n) = 0$ for more negative $s$], L25 gives a complete proof in the following cases: (i) any $s < 0$ and $n + 2s \equiv 0 \bmod 4$ or $n + 2s \equiv 2 \bmod 4$; (ii) any $s < -1$ and $n + 2s \equiv 1 \bmod 4$; and (iii) any $s < -2$ and $n + 2s \equiv 3 \bmod 4$. For the exceptional cases where $s = -1$ and $n + 2s \equiv 1 \bmod 4$ or $n + 2s \equiv 3 \bmod 4$, or where $s = -2$ and $n + 2s \equiv 3 \bmod 4$, see Properties 6.1 and 6.2 of L25. Once established, the foregoing additional properties would complete the formal proof of the global unimodality of $H_s(n)$ at $s=0$ and would confirm the blimpy shape of the log number of heady-$s$ bit strings of any length $n$.

**Appendix** A.

We show $f_{a,b}(x) = 1/\{(a+b/x)^2 - 1\}$ is convex over the open interval $0 < x < b/(1-a)$ for any $0 \le a < 1$ and $b > 0$. The first two derivatives are

$$f'_{a,b}(x) = -\{(a+b/x)^2 - 1\}^{-2} \cdot 2(a+b/x)(-b/x^2) \;=\; 2b(ax^{-2} + bx^{-3}) f_{a,b}(x)^2 \;>\; 0$$

and

$$(2b)^{-1} f''_{a,b}(x) = (ax^{-2} + bx^{-3}) \cdot 2 f_{a,b}(x) \{2b(ax^{-2} + bx^{-3}) f_{a,b}(x)^2\} \;-\; f_{a,b}(x)^2 (2ax^{-3} + 3bx^{-4}) \,.$$

The latter expression is positive if and only if $4b(ax^{-2} + bx^{-3})^2 f_{a,b}(x) \;>\; (2ax^{-3} + 3bx^{-4})$. Multiplying both sides by $(a+b/x)^2 - 1$ and then by $x^6$ on the left and $x^4 \cdot x^2$ on the right, the inequality holds if and only if

$$4b(ax+b)^2 \;>\; (2ax+3b)\{(ax+b)^2 - x^2\}$$

if and only if

$$(ax+b)^2(b - 2ax) + x^2(2ax + 3b) \;>\; 0$$

if and only if

$$a^2bx^2 + 2ab^2x + b^3 - 2a^3x^3 - 4a^2bx^2 - 2ab^2x + 2ax^3 + 3bx^2 \;>\; 0$$

if and only if

$$2a(1-a^2)x^3 + (a^2b - 4a^2b + 3b)x^2 + b^3 \;=\; 2a(1-a^2)x^3 + 3b(1-a^2)x^2 + b^3 \;>\; 0$$

which is true for $0 \le a < 1$, $b > 0$, and $x > 0$. □

**Appendix B.**

We prove Lemma 2.1. To show $\Delta\{s_j t_j\} = s_{j+1}t_{j+1} - s_j t_j > 0$, we have

$$s_{j+1}t_{j+1} - s_j t_j = s_{j+1}t_{j+1} - s_{j+1}t_j + s_{j+1}t_j - s_j t_j = s_{j+1}\Delta t_j + (\Delta s_j) t_j > 0$$

since $\{s_j\}$ and $\{t_j\}$ are positive increasing. To show $\Delta\{s_j t_j\} > \Delta\{s_{j-1}t_{j-1}\}$, the forward difference on the left is as above, and that on the right is $s_j \Delta t_{j-1} + (\Delta s_{j-1}) t_{j-1}$. But $s_{j+1}\Delta t_j > s_j \Delta t_j > s_j \Delta t_{j-1}$ since $\Delta t_j > 0$, $s_j$ increases, and $\Delta t_j$ increases. Similarly, $(\Delta s_j) t_j > (\Delta s_j) t_{j-1} > (\Delta s_{j-1}) \Delta t_{j-1}$ since $\Delta s_j > 0$, $t_j$ increases, and $\Delta s_j$ increases. Hence $\Delta^2\{s_j t_j\} = \Delta\Delta\{s_{j-1}t_{j-1}\} = \Delta\{s_j t_j\} - \Delta\{s_{j-1}t_{j-1}\} > 0$. □

**Appendix C.**

We prove Lemma 4.1. From (4.2) we have

$$\begin{aligned}\log(1+LHS) &= \log\prod_{i=0}^{2a}\frac{(4a^2+a-1-i)}{(4a^2-a-3-i)}+\log\prod_{i=0}^{2a+1}\frac{(4a^2+a-i)}{(4a^2+3a+2-i)}\\ &=\sum_{i=0}^{2a}\log\frac{1+(\alpha/4)\{1-(i+1)\alpha\}}{1-(\alpha/4)\{1+(i+3)\alpha\}}+\sum_{i=0}^{2a}\log\frac{1+(\alpha/4)(1-i\alpha)}{1+(\alpha/4)\{3-(i-2)\alpha\}}+t_{2a+1}\end{aligned} \tag{C.1}$$

where the final term $t_{2a+1}=\log\dfrac{1+(\alpha/4)\{1-(2a+1)\alpha\}}{1+(\alpha/4)\{3-(2a-1)\alpha\}}=\log\dfrac{1-(\alpha/4)(1+\alpha)}{1+(\alpha/4)(1+\alpha)}$ has expansion

$$\begin{aligned}t_{2a+1} &= -\sum_{\nu=1}^{\infty}\{1+(-1)^{\nu-1}\}(\alpha/4)^{\nu}(1+\alpha)^{\nu}/\nu=-2\{(\alpha/4)(1+\alpha)+(\alpha/4)^3(1+\alpha)^3/3\}+O(\alpha^5)\\ &=-\{(\alpha/2)+(\alpha^2/2)+(\alpha^3/96)\}+o(\alpha^3)\,.\end{aligned} \tag{C.2}$$

For $i=0,\ldots,2a$, the individual summands preceding $t_{2a+1}$ in (C.1) have asymptotic expansions

$$\begin{aligned}&(\frac{\alpha}{4})\{1-(i+1)\alpha\} - (\frac{\alpha}{4})^2\{1-(i+1)\alpha\}^2/2 + (\frac{\alpha}{4})^3\{1-(i+1)\alpha\}^3/3 - (\frac{\alpha}{4})^4\{1-(i+1)\alpha\}^4/4 +-\cdots\\ &+ (\frac{\alpha}{4})\{1+(i+3)\alpha\} + (\frac{\alpha}{4})^2\{1+(i+3)\alpha\}^2/2 + (\frac{\alpha}{4})^3\{1+(i+3)\alpha\}^3/3 + (\frac{\alpha}{4})^4\{1+(i+3)\alpha\}^4/4+\cdots\\ &+ (\frac{\alpha}{4})(1-i\alpha) - (\frac{\alpha}{4})^2(1-i\alpha)^2/2 + (\frac{\alpha}{4})^3(1-i\alpha)^3/3 - (\frac{\alpha}{4})^4(1-i\alpha)^4/4+-\cdots\\ &- (\frac{\alpha}{4})\{3-(i-2)\alpha\} + (\frac{\alpha}{4})^2\{3-(i-2)\alpha\}^2/2 - (\frac{\alpha}{4})^3\{3-(i-2)\alpha\}^3/3 + (\frac{\alpha}{4})^4\{3-(i-2)\alpha\}^4/4 -+\cdots.\end{aligned}$$

The four terms beginning with $(\frac{\alpha}{4})$ sum to zero for each $i$. The sum of the four terms beginning with $(\frac{\alpha}{4})^2$, after expanding the squared braces, is

$$\begin{aligned}\frac{1}{2}(\frac{\alpha}{4})^2\{&-1 \quad +2(i+1)\alpha \quad -(i+1)^2\alpha^2\\ &+1 \quad +2(i+3)\alpha \quad +(i+3)^2\alpha^2\\ &-1 \quad +2i\alpha \quad -i^2\alpha^2\\ &+9 \quad -6(i-2)\alpha \quad +(i-2)^2\alpha^2\}\\ =\ \frac{1}{2}(\frac{\alpha}{4})^2&(8+20\alpha+12\alpha^2)\ =\ \frac{\alpha^2}{4}+\frac{5\alpha^3}{8}+\frac{3\alpha^4}{8}+o(\alpha^4)\end{aligned}$$

with no dependence on $i$. Then, summing over $i = 0,...,2a$ as in (C.1) yields

$$\{\frac{\alpha^2}{4} + \frac{5\alpha^3}{8} + \frac{3\alpha^4}{8} + o(\alpha^4)\}\,(2a+1) = \frac{\alpha}{2} + \frac{\alpha^2}{4} + \frac{5\alpha^2}{4} + \frac{5\alpha^3}{8} + \frac{3\alpha^3}{4} + o(\alpha^3)$$
$$= \frac{\alpha}{2} + \frac{3\alpha^2}{2} + \frac{11\alpha^3}{8} + o(\alpha^3)\,. \tag{C.3}$$

The sum of the four terms beginning with $(\frac{\alpha}{4})^3$, after expanding the cubed braces, is

$$\begin{aligned}\frac{1}{3}(\frac{\alpha}{4})^3\{1 &\quad -3(i+1)\alpha \quad +3(i+1)^2\alpha^2 \quad -(i+1)^3\alpha^3\\ +1 &\quad +3(i+3)\alpha \quad +3(i+3)^2\alpha^2 \quad +(i+3)^3\alpha^3\\ +1 &\quad -3i\alpha \quad +3i^2\alpha^2 \quad -i^3\alpha^3\\ -27 &\quad +27(i-2)\alpha \quad -9(i-2)^2\alpha^2 \quad +(i-2)^3\alpha^3\,\}.\end{aligned}$$

Expanding the powers involving $i$ and retaining in the second line below only the terms from the first line that become of order $O(\alpha^3)$ or less after summing over $i = 0,...,2a$, we obtain

$$\frac{1}{3}(\frac{\alpha}{4})^3\{-24 + 24i\alpha - 48\alpha + 60i\alpha^2 - 6\alpha^2 + 36i\alpha^3 + 18\alpha^3\}$$
$$= -\frac{\alpha^3}{8} + \frac{i\alpha^4}{8} - \frac{\alpha^4}{4} + \frac{5i\alpha^5}{16} \;+\; \text{higher order terms.}$$

Then, summing over $i = 0,...,2a$, yields

$$\{-\frac{\alpha^3}{8} - \frac{\alpha^4}{4}\}\,(2a+1) \;+\; \{\frac{\alpha^4}{8} + \frac{5\alpha^5}{16}\}(2a)(2a+1)/2 \;+\; o(\alpha^3)$$
$$= -\frac{\alpha^2}{4} - \frac{\alpha^3}{8} - \frac{\alpha^3}{2} \;+\; \frac{\alpha^2}{4} + \frac{\alpha^3}{8} + \frac{5\alpha^3}{8} \;+\; o(\alpha^3) = \frac{\alpha^3}{8} + o(\alpha^3). \tag{C.4}$$

The sum of the four terms beginning with $(\frac{\alpha}{4})^4$, after expanding the fourth-power braces, is

$$\begin{aligned}\frac{1}{4}(\frac{\alpha}{4})^4\{-1 &\quad +4(i+1)\alpha \quad -6(i+1)^2\alpha^2 \quad +4(i+1)^3\alpha^3 \quad -(i+1)^4\alpha^4\\ +1 &\quad +4(i+3)\alpha \quad +6(i+3)^2\alpha^2 \quad +4(i+3)^3\alpha^3 \quad +(i+3)^4\alpha^4\\ -1 &\quad +4i\alpha \quad -6i^2\alpha^2 \quad +4i^3\alpha^3 \quad -i^4\alpha^4\\ +81 &\quad -4\cdot 27(i-2)\alpha \quad +6\cdot 9(i-2)^2\alpha^2 \quad -4\cdot 3(i-2)^3\alpha^3 \quad +(i-2)^4\alpha^4\,\}.\end{aligned}$$

Expanding the powers involving $i$ and retaining in the second line below only the terms from the first line that become of order $O(\alpha^3)$ or less after summing over $i = 0,...,2a$, we obtain

$$\frac{1}{4}(\frac{\alpha}{4})^4\{80+(-96i+232)\alpha+(48i^2-192i+264)\alpha^2+(120i^2-24i+208)\alpha^3+(72i^2+72i+96)\alpha^4\}$$

$$=\frac{5\alpha^4}{64}-\frac{3i\alpha^5}{32}+\frac{3i^2\alpha^6}{64}+ \text{ higher order terms.}$$

Then, summing over $i = 0,...,2a$, yields

$$(\frac{5\alpha^4}{64})\,(2a+1) \;-\; (\frac{3\alpha^5}{32})(2a)(2a+1)/2 \;+\; (\frac{3\alpha^6}{64})(2a)(2a+1)(4a+1)/6 \;+\; o(\alpha^3)$$
$$=\frac{5\alpha^3}{32}-\frac{3\alpha^3}{16}+\frac{\alpha^3}{8} \;+\; o(\alpha^3)=\frac{3\alpha^3}{32} \;+\; o(\alpha^3). \tag{C.5}$$

It is easy to check that all higher powers in the logarithmic expansions only yield terms of order $o(\alpha^3)$. Therefore, adding (C.2), (C.3), (C.4), and (C.5) yields

$$\log(1+LHS)=-\left(\frac{\alpha}{2}+\frac{\alpha^2}{2}+\frac{\alpha^3}{96}\right) \;+\; \left(\frac{\alpha}{2}+\frac{3\alpha^2}{2}+\frac{11\alpha^3}{8}\right) \;+\; \left(\frac{\alpha^3}{8}\right) \;+\; \left(\frac{3\alpha^3}{32}\right) \;+\; o(\alpha^3)$$
$$=\alpha^2+(19/12)\alpha^3 \;+\; o(\alpha^3)$$

as was to be shown. Finally, $\exp(u)-1=u+O(u^2)$ as $u \to 0$, so with $u=\log(1+LHS)$,

$$LHS=\exp\{\log(1+LHS)\}-1=\log(1+LHS)+O[\{\log(1+LHS)\}^2]$$
$$=\alpha^2+(19/12)\alpha^3+o(\alpha^3)+O[\{\alpha^2+(19/12)\alpha^3+o(\alpha^3)\}^2]$$
$$=\alpha^2+(19/12)\alpha^3+o(\alpha^3) \qquad \text{as well.} \quad \square$$

**Appendix D.**

We prove Lemma 4.2 by showing that the factor of *RHS* in square brackets at (4.3) equals $\alpha^3/2+o(\alpha^3)$ and that the remaining product-ratio of binomial coefficients converges to Euler's constant *e* as $\alpha \to 0$. Thus $RHS=\{\alpha^3/2+o(\alpha^3)\}\{e+o(1)\}=(e/2)\alpha^3+o(\alpha^3)$.

For the factor of (4.3) in square brackets, dividing numerator and denominator by $64a^7$ yields

$$\frac{\{32a^2(a^2-1)-6\}\cdot(4a+1)/(4a)}{(4a^2+2)(4a^2+1)(4a^2-3)\cdot(a+1)}=\left(\frac{\alpha^3}{2}\right)\left(\frac{(1-\alpha^2-\frac{3\alpha^4}{16})(1+\frac{\alpha}{4})}{(1+\frac{\alpha^2}{2})(1+\frac{\alpha^2}{4})(1-\frac{3\alpha^2}{4})(1+\alpha)}\right)$$

$$=\left(\frac{\alpha^3}{2}\right)\left(\frac{1+\frac{\alpha}{4}+o(\alpha)}{1+\alpha+o(\alpha)}\right)=\frac{\alpha^3}{2}+o(\alpha^3).$$

For the remaining binomial coefficients in (4.3), we could apply Stirling's approximation to each one, but it will be more expedient here to do some preliminary manipulation of the whole product-ratio factor of binomial coefficients, leading to (D.1) and (D.2) below. The product-ratio factor is

$$\frac{\binom{8a^2-1}{4a^2-1}\binom{8a^2-5}{4a^2-1}}{\binom{8a^2+2a-1}{4a^2+a-1}\binom{8a^2-2a-5}{4a^2+a-1}}=\frac{\binom{8a^2-1}{4a^2}\binom{8a^2-5}{4a^2-4}}{\binom{8a^2+2a-1}{4a^2+a}\binom{8a^2-2a-5}{4a^2-3a-4}}$$

$$=\frac{\dfrac{(8a^2-1)^{(4a^2)}(8a^2-5)^{(4a^2-4)}}{(4a^2)^{(4a^2)}(4a^2-4)^{(4a^2-4)}}}{\dfrac{(8a^2+2a-1)^{(4a^2+a)}(8a^2-2a-5)^{(4a^2-3a-4)}}{(4a^2+a)^{(4a^2+a)}(4a^2-3a-4)^{(4a^2-3a-4)}}}=\frac{\prod_{i=0}^{4a^2-1}\frac{8a^2-1-i}{4a^2-i}\prod_{i=0}^{4a^2-5}\frac{8a^2-5-i}{4a^2-4-i}}{\prod_{i=0}^{4a^2+a-1}\frac{8a^2+2a-1-i}{4a^2+a-i}\prod_{i=0}^{4a^2-3a-5}\frac{8a^2-2a-5-i}{4a^2-3a-4-i}}.$$

Changing the index of summation to $i'=4a^2-i$ in the first numerator factor, $i'=4a^2-4-i$ in the second numerator factor, $i'=4a^2+a-i$ in the first denominator factor, and $i'=4a^2-3a-4-i$ in the second denominator factor, then dropping the primes, the expression becomes

$$\frac{\prod_{i=1}^{4a^2}\left\{1+\frac{4a^2-1}{i}\right\}\prod_{i=1}^{4a^2-4}\left\{1+\frac{4a^2-1}{i}\right\}}{\prod_{i=1}^{4a^2+a}\left\{1+\frac{4a^2+a-1}{i}\right\}\prod_{i=1}^{4a^2-3a-4}\left\{1+\frac{4a^2+a-1}{i}\right\}}$$

$$=\frac{\prod_{i=1}^{4a^2-4}\left\{1+\frac{4a^2-1}{i}\right\}\prod_{i=4a^2-3}^{4a^2}\left\{1+\frac{4a^2-1}{i}\right\}\prod_{i=1}^{4a^2-4}\left\{1+\frac{4a^2-1}{i}\right\}}{\prod_{i=1}^{4a^2-3a-4}\left\{1+\frac{4a^2+a-1}{i}\right\}\prod_{i=4a^2-3a-3}^{4a^2+a}\left\{1+\frac{4a^2+a-1}{i}\right\}\prod_{i=1}^{4a^2-3a-4}\left\{1+\frac{4a^2+a-1}{i}\right\}}$$

$$=\frac{\prod_{i=0}^{3}\left\{1+\frac{4a^2-1}{4a^2-3+i}\right\}\prod_{i=1}^{4a^2-4}\left\{1+\frac{4a^2-1}{i}\right\}^2}{\prod_{i=0}^{4a+3}\left\{1+\frac{4a^2+a-1}{4a^2-3a-3+i}\right\}\prod_{i=1}^{4a^2-3a-4}\left\{1+\frac{4a^2+a-1}{i}\right\}^2}$$

$$=\frac{\prod_{i=0}^{3}\left\{1+\frac{4a^2-1}{4a^2-3+i}\right\}\prod_{i=1}^{4a^2-3a-4}\left\{1+\frac{4a^2-1}{i}\right\}^2\prod_{i=4a^2-3a-3}^{4a^2-4}\left\{1+\frac{4a^2-1}{i}\right\}^2}{\prod_{i=0}^{3a-1}\left\{1+\frac{4a^2+a-1}{4a^2-3a-3+i}\right\}\prod_{i=3a}^{4a+3}\left\{1+\frac{4a^2+a-1}{4a^2-3a-3+i}\right\}\prod_{i=1}^{4a^2-3a-4}\left\{1+\frac{4a^2+a-1}{i}\right\}^2}$$

$$=\frac{\prod_{i=0}^{3}\left\{1+\frac{4a^2-1}{4a^2-3+i}\right\}\prod_{i=0}^{3a-1}\left\{1+\frac{4a^2-1}{4a^2-3a-3+i}\right\}^2\prod_{i=1}^{4a^2-3a-4}\left\{\frac{1+\frac{4a^2-1}{i}}{1+\frac{4a^2+a-1}{i}}\right\}^2}{\prod_{i=0}^{3a-1}\left\{1+\frac{4a^2+a-1}{4a^2-3a-3+i}\right\}\prod_{i=0}^{a+3}\left\{1+\frac{4a^2+a-1}{4a^2-3+i}\right\}}$$

$$=\frac{\prod_{i=0}^{3}\left\{1+\frac{1-\alpha^2/4}{1+(\alpha^2/4)(i-3)}\right\}\prod_{i=0}^{3a-1}\left\{1+\frac{1-\alpha^2/4}{1-(3\alpha/4)+(\alpha^2/4)(i-3)}\right\}^2\prod_{i=1}^{4a^2-3a-4}\left\{\frac{1+\frac{1-\alpha^2/4}{(\alpha^2/4)i}}{1+\frac{1+(\alpha/4)-\alpha^2/4}{(\alpha^2/4)i}}\right\}^2}{\prod_{i=0}^{3a-1}\left\{1+\frac{1+(\alpha/4)-\alpha^2/4}{1-(3\alpha/4)+(\alpha^2/4)(i-3)}\right\}\prod_{i=0}^{a+3}\left\{1+\frac{1+(\alpha/4)-\alpha^2/4}{1+(\alpha^2/4)(i-3)}\right\}}.$$

Now the third factor in the above numerator is

$$\prod_{i=1}^{4a^2-3a-4}\left\{\frac{1+\dfrac{1-\alpha^2/4}{(\alpha^2/4)i}}{1+\dfrac{1+(\alpha/4)-\alpha^2/4}{(\alpha^2/4)i}}\right\}^2 = \prod_{i=1}^{4a^2-3a-4}\left\{\frac{1+(\alpha^2/4)(i-1)}{1+(\alpha/4)+(\alpha^2/4)(i-1)}\right\}^2$$

$$= \prod_{i=1}^{4a^2-3a-4}\left\{1-\frac{\alpha/4}{1+(\alpha/4)+(\alpha^2/4)(i-1)}\right\}^2$$

and the second factor in the denominator is

$$\prod_{i=0}^{a+3}\left\{1+\frac{1+(\alpha/4)-\alpha^2/4}{1+(\alpha^2/4)(i-3)}\right\} = \prod_{i=0}^{3}\left\{1+\frac{1+(\alpha/4)-\alpha^2/4}{1+(\alpha^2/4)(i-3)}\right\}\prod_{i=4}^{a+3}\left\{1+\frac{1+(\alpha/4)-\alpha^2/4}{1+(\alpha^2/4)(i-3)}\right\}$$

$$= \prod_{i=0}^{3}\left\{1+\frac{1+(\alpha/4)-\alpha^2/4}{1+(\alpha^2/4)(i-3)}\right\}\prod_{i=1}^{a}\left\{1+\frac{1+(\alpha/4)-\alpha^2/4}{1+(\alpha^2/4)i}\right\}.$$

Substituting, then simplifying the leading $\Pi_0^3$ factor and shifting the index of the two $\Pi_0^{3a-1}$ factors, the product-ratio expression equals

$$\frac{\displaystyle\prod_{i=0}^{3}\left\{\frac{1+\dfrac{1-\alpha^2/4}{1+(\alpha^2/4)(i-3)}}{1+\dfrac{1+(\alpha/4)-\alpha^2/4}{1+(\alpha^2/4)(i-3)}}\right\}\prod_{i=0}^{3a-1}\left\{1+\frac{1-\alpha^2/4}{1-(3\alpha/4)+(\alpha^2/4)(i-3)}\right\}^2\prod_{i=1}^{4a^2-3a-4}\left\{1-\frac{\alpha/4}{1+(\alpha/4)+(\alpha^2/4)(i-1)}\right\}^2}{\displaystyle\prod_{i=1}^{a}\left\{1+\frac{1+(\alpha/4)-\alpha^2/4}{1+(\alpha^2/4)i}\right\}\prod_{i=0}^{3a-1}\left\{1+\frac{1+(\alpha/4)-\alpha^2/4}{1-(3\alpha/4)+(\alpha^2/4)(i-3)}\right\}}$$

$$= \prod_{i=0}^{3}\left\{1-\frac{\alpha/4}{2+(\alpha/4)+(\alpha^2/4)(i-4)}\right\}$$

$$\times\frac{\displaystyle\prod_{i=1}^{3a}\left\{1+\frac{1-\alpha^2/4}{1-(3\alpha/4)+(\alpha^2/4)(i-4)}\right\}^2\prod_{i=1}^{4a^2-3a-4}\left\{1-\frac{\alpha/4}{1+(\alpha/4)+(\alpha^2/4)(i-1)}\right\}^2}{\displaystyle\prod_{i=1}^{a}\left\{1+\frac{1+(\alpha/4)-\alpha^2/4}{1+(\alpha^2/4)i}\right\}\prod_{i=1}^{3a}\left\{1+\frac{1+(\alpha/4)-\alpha^2/4}{1-(3\alpha/4)+(\alpha^2/4)(i-4)}\right\}}. \tag{D.1}$$

Thus, the log-product-ratio is of the form

$$\log\frac{\dbinom{8a^2-1}{4a^2-1}\dbinom{8a^2-5}{4a^2-1}}{\dbinom{8a^2+2a-1}{4a^2+a-1}\dbinom{8a^2-2a-5}{4a^2+a-1}} = L_0 - L_1 + 2L_2 - L_3 + 2L_4\,, \tag{D.2}$$

where

$$L_0 := \sum_{i=0}^{3} \log\left\{1 - \frac{\alpha/4}{2 + (\alpha/4) + (\alpha^2/4)(i-4)}\right\},$$

$$L_1 := \sum_{i=1}^{a} \log\left\{1 + \frac{1 + (\alpha/4) - (\alpha^2/4)}{1 + (\alpha^2/4)i}\right\},$$

$$L_2 := \sum_{i=1}^{3a} \log\left\{1 + \frac{1 - (\alpha^2/4)}{1 - (3\alpha/4) + (\alpha^2/4)(i-4)}\right\},$$

$$L_3 := \sum_{i=1}^{3a} \log\left\{1 + \frac{1 + (\alpha/4) - (\alpha^2/4)}{1 - (3\alpha/4) + (\alpha^2/4)(i-4)}\right\},$$

and

$$L_4 := \sum_{i=1}^{4a^2-3a-4} \log\left\{1 - \frac{\alpha/4}{1 + (\alpha/4) + (\alpha^2/4)(i-1)}\right\}.$$

With re-expression (D.2) in hand, below we'll only need a Stirling-type approximation for $L_4$.

Claim D.1. As $a \to \infty$ and $\alpha \to 0$, (0) $L_0 = o(1)$;

(1) $L_1 = (\log 2)a + (1/16) + o(1)$;

(2) $L_2 = (3\log 2)a + (9/16) + o(1)$;

(3) $L_3 = (3\log 2)a + (15/16) + o(1)$;

and (4) $L_4 = (-\log 2)a + (7/16) + o(1)$.

Thus $L_0 - L_1 + 2L_2 - L_3 + 2L_4 = (0-1+6-3-2)(\log 2)a + (0-1+18-15+14)/16 + o(1) = 1 + o(1)$. Hence the product-ratio converges to *e*, as stated in Lemma 4.2.

Proof of Claim D.1. Claim (0) is obvious. For (1), write

$$L_1 = \sum_{i=1}^{a} \log\left\{1 + \frac{1 + (\alpha/4) - (\alpha^2/4)}{1 + (\alpha^2/4)i}\right\} = \sum_{i=1}^{a} \log\left\{\frac{2 + (\alpha/4)(1-\alpha) + (\alpha^2/4)i}{1 + (\alpha^2/4)i}\right\}$$

$$= \sum_{i=1}^{a} \log\{2 + (\alpha/4)(1-\alpha) + (\alpha^2/4)i\} - \sum_{i=1}^{a} \log\{1 + (\alpha^2/4)i\}$$

$$= (\log 2)a + \sum_{i=1}^{a} \log\{1 + (\alpha^2/8)i + (\alpha/8)(1-\alpha)\} - \sum_{i=1}^{a} \log\{1 + (\alpha^2/4)i\}$$

$$= (\log 2)a + \sum_{i=1}^{a} \left[\begin{array}{c} \{(\alpha^2/8)i + (\alpha/8)(1-\alpha)\} - \{(\alpha^2/8)i + (\alpha/8)(1-\alpha)\}^2/2 + - \cdots \\ - \{(\alpha^2/4)i\} + \{(\alpha^2/4)i\}^2/2 - + \cdots \end{array}\right]$$

$$= (\log 2)a - (\alpha^2/8)a(a+1)/2 + (\alpha/8)(1-\alpha)a + o(1)$$

$$= (\log 2)a - (1/16) + (1/8) + o(1) \;=\; (\log 2)a + (1/16) + o(1).$$

For (2), write

$$L_2 = \sum_{i=1}^{3a} \log\left\{1 + \frac{1-(\alpha^2/4)}{1-(3\alpha/4)+(\alpha^2/4)(i-4)}\right\} = \sum_{i=1}^{3a} \log\left\{1 + \frac{1-(\alpha^2/4)}{1-(3\alpha/4)-\alpha^2+(\alpha^2/4)i)}\right\}$$

$$= \sum_{i=1}^{3a} \log\left\{\frac{2-(3\alpha/4)-(5\alpha^2/4)+(\alpha^2/4)i}{1-(3\alpha/4)-\alpha^2+(\alpha^2/4)i)}\right\}$$

$$= (3\log 2)a + \sum_{i=1}^{3a} \log\{1+(\alpha^2/8)i-(3\alpha/8)-(5\alpha^2/8)\} - \sum_{i=1}^{3a} \log\{1+(\alpha^2/4)i-(3\alpha/4)-\alpha^2\}$$

$$= (3\log 2)a + \sum_{i=1}^{3a} \left[\{(\alpha^2/8)i-(3\alpha/8)-(5\alpha^2/8)\} - + \cdots \quad - \{(\alpha^2/4)i-(3\alpha/4)-\alpha^2\} + - \cdots\right]$$

$$= (3\log 2)a + (\alpha^2/8)(3a)(3a+1)/2 - (3\alpha/8)(3a) - (\alpha^2/4)(3a)(3a+1)/2 + (3\alpha/4)(3a) + o(1)$$

$$= (3\log 2)a + (9/16) - (9/8) - (9/8) + (9/4) + o(1) \ = \ (3\log 2)a + (9/16) + o(1)\,.$$

For (3), which is similar to (2), write

$$L_3 = \sum_{i=1}^{3a} \log\left\{1 + \frac{1+(\alpha/4)-(\alpha^2/4)}{1-(3\alpha/4)+(\alpha^2/4)(i-4)}\right\} = \sum_{i=1}^{3a} \log\left\{1 + \frac{1+(\alpha/4)-(\alpha^2/4)}{1-(3\alpha/4)-\alpha^2+(\alpha^2/4)i)}\right\}$$

$$= \sum_{i=1}^{3a} \log\left\{\frac{2-(\alpha/2)-(5\alpha^2/4)+(\alpha^2/4)i}{1-(3\alpha/4)-\alpha^2+(\alpha^2/4)i)}\right\}$$

$$= (3\log 2)a + \sum_{i=1}^{3a} \log\{1+(\alpha^2/8)i-(\alpha/4)-(5\alpha^2/8)\} - \sum_{i=1}^{3a} \log\{1+(\alpha^2/4)i-(3\alpha/4)-\alpha^2\}$$

$$= (3\log 2)a + \sum_{i=1}^{3a} \left[\{(\alpha^2/8)i-(\alpha/4)-(5\alpha^2/8)\} - + \cdots \quad - \{(\alpha^2/4)i-(3\alpha/4)-\alpha^2\} + - \cdots\right]$$

$$= (3\log 2)a + (\alpha^2/8)(3a)(3a+1)/2 - (\alpha/4)(3a) - (\alpha^2/4)(3a)(3a+1)/2 + (3\alpha/4)(3a) + o(1)$$

$$= (3\log 2)a + (9/16) - (3/4) - (9/8) + (9/4) + o(1) \ = \ (3\log 2)a + (15/16) + o(1)\,.$$

For (4), a similar approach becomes more arduous because $L_4$ has $4a^2-3a-4$ summands and infinitely many terms in the logarithmic expansion contribute to the zero-order term of $L_4$. Instead, we use a Stirling-like integral approximation for the sum with ½ continuity correction. Write

$$L_4 = \sum_{i=1}^{4a^2-3a-4} \log\left\{1 - \frac{\alpha/4}{1+(\alpha/4)+(\alpha^2/4)(i-1)}\right\} = \sum_{i=1}^{4a^2-3a-4} \log\left\{1 - \frac{\alpha/4}{1+(\alpha/4)(1-\alpha)+(\alpha^2/4)i}\right\}$$

$$= \int_{1/2}^{4a^2-3a-4+1/2} \log\{1 - \frac{\gamma_0}{1+\gamma_1+\gamma_2 x}\} dx + o(1)\,,$$

where $\gamma_0 = \alpha/4$, $\gamma_1 = (\alpha/4)(1-\alpha)$, and $\gamma_2 = \alpha^2/4$. Now change variables to $y = 1+\gamma_1+\gamma_2 x$, so that the integral equals $I = \int_{1/2}^{4a^2-3a-4+1/2} \log\{1 - \frac{\gamma_0}{1+\gamma_1+\gamma_2 x}\}dx = \gamma_2^{-1}\int_{y_L}^{y_U} \log\{1-\frac{\gamma_0}{y}\}dy$, where

$$
\begin{aligned}
y_L &= 1+\gamma_1+\gamma_2/2 \\
&= 1+(\alpha/4)-(\alpha^2/4)+(\alpha^2/8) \\
&= 1+(\alpha/4)-(\alpha^2/8)
\end{aligned}
$$

and

$$
\begin{aligned}
y_U &= 1+\gamma_1+\gamma_2(4a^2-3a-4+½) \\
&= 1+(\alpha/4)-(\alpha^2/4)+(\alpha^2/4)(4a^2-3a-4+½) \\
&= 1+(\alpha/4)-(\alpha^2/4)+1-(3\alpha/4)-\alpha^2+(\alpha^2/8) \\
&= 2-(\alpha/2)-(9\alpha^2/8).
\end{aligned}
$$

The indefinite integral of $\log\{1-\frac{\gamma_0}{y}\}$ is $y\log(1-\gamma_0/y)+\gamma_0\{1-\log(y-\gamma_0)\}$, so

$$\gamma_2 I = y_U \log(1-\gamma_0/y_U)+\gamma_0\{1-\log(y_U-\gamma_0)\} \;-\; y_L\log(1-\gamma_0/y_L)-\gamma_0\{1-\log(y_L-\gamma_0)\}.$$

Now for either $y = y_U$ or $y_L$,

$$y\log(1-\gamma_0/y) = -\gamma_0 - \frac{\gamma_0^2}{2y}+o(\gamma_0^2) = -(\alpha/4)-(\frac{\alpha^2}{32y})+o(\alpha^2)$$

and

$$\gamma_0\log(y-\gamma_0) = \gamma_0\log y+\gamma_0\log(1-\gamma_0/y) = \gamma_0\log y-\gamma_0^2/y+o(\gamma_0^2) = (\alpha/4)\log y-(\frac{\alpha^2}{16y})+o(\alpha^2).$$

Therefore,

$$
\begin{aligned}
\gamma_2 I &= \left[-\frac{\alpha}{4}-\frac{\alpha^2}{32y_U}+\frac{\alpha}{4}-\frac{\alpha}{4}\log y_U+\frac{\alpha^2}{16y_U}+o(\alpha^2)\right] - \left[-\frac{\alpha}{4}-\frac{\alpha^2}{32y_L}+\frac{\alpha}{4}-\frac{\alpha}{4}\log y_L+\frac{\alpha^2}{16y_L}+o(\alpha^2)\right] \\
&= \frac{\alpha^2}{32y_U}-\frac{\alpha}{4}\log y_U-\frac{\alpha^2}{32y_L}+\frac{\alpha}{4}\log y_L+o(\alpha^2).
\end{aligned}
$$

But

$$\frac{\alpha^2}{32y_U} = \frac{\alpha^2}{32\{2-\frac{\alpha}{2}-\frac{9\alpha^2}{8}\}} = \frac{\alpha^2}{64}+o(\alpha^2),$$

$$-\frac{\alpha}{4}\log y_U = -\frac{\alpha}{4}\log\{2-\frac{\alpha}{2}-\frac{9\alpha^2}{8}\} = -\frac{\alpha}{4}\log 2-\frac{\alpha}{4}\log\{1-\frac{\alpha}{4}-\frac{9\alpha^2}{16}\} = -\frac{\alpha}{4}\log 2+\frac{\alpha^2}{16}+o(\alpha^2),$$

$$-\frac{\alpha^2}{32y_L} = -\frac{\alpha^2}{32\{1+\frac{\alpha}{4}-\frac{\alpha^2}{8}\}} = -\frac{\alpha^2}{32}+o(\alpha^2)\,,$$

and

$$\frac{\alpha}{4}\log y_L = \frac{\alpha}{4}\log\{1+\frac{\alpha}{4}-\frac{\alpha^2}{8}\} = \frac{\alpha}{4}\{\frac{\alpha}{4}-\frac{\alpha^2}{8}+o(\alpha^2)\} = \frac{\alpha^2}{16}+o(\alpha^2)\,.$$

Thus,

$$\gamma_2 I = \frac{\alpha^2}{64}-\frac{\alpha}{4}\log 2+\frac{\alpha^2}{16}-\frac{\alpha^2}{32}+\frac{\alpha^2}{16}+o(\alpha^2) = -\frac{\alpha}{4}\log 2+(\frac{1}{64}+\frac{1}{16}-\frac{1}{32}+\frac{1}{16})\alpha^2+o(\alpha^2)$$

$$= -\frac{\alpha}{4}\log 2+\frac{7\alpha^2}{64}+o(\alpha^2)\,.$$

Therefore,

$$L_4 = I+o(1) = 4a^2\left\{-\frac{\alpha}{4}\log 2+\frac{7\alpha^2}{64}+o(\alpha^2)\right\}+o(1) = (-\log 2)a+(7/16)+o(1)\,.$$

This concludes the proof of Claim D.1 and Lemma 4.2. □

Remark D.1. The convergence of the product-ratio of binomial coefficients to $e$ as $\alpha \to 0$ is quite slow. For example, for (shifted) $a$=10, 100, 1000, and 5000, the values of

$$\frac{\binom{8a^2-1}{4a^2-1}\binom{8a^2-5}{4a^2-1}}{\binom{8a^2+2a-1}{4a^2+a-1}\binom{8a^2-2a-5}{4a^2+a-1}}$$ equal 3.2737, 2.7665, 2.7230, and 2.7189, respectively.

**Acknowledgment.**

I thank Professor Cheng-Shiun Leu for helpful discussions and critical reading during the preparation of this report.

**References.**

Figure 4.1

Illustrating the pairwise inequalities for $d(j)Q(j)$ in the case $n=1024$, $l=256$, $u=0$, $a=8$.

Red boxes are $d(j+a+1)Q(j+a+1)$. Blue diamonds are $-d(j)Q(j)$.

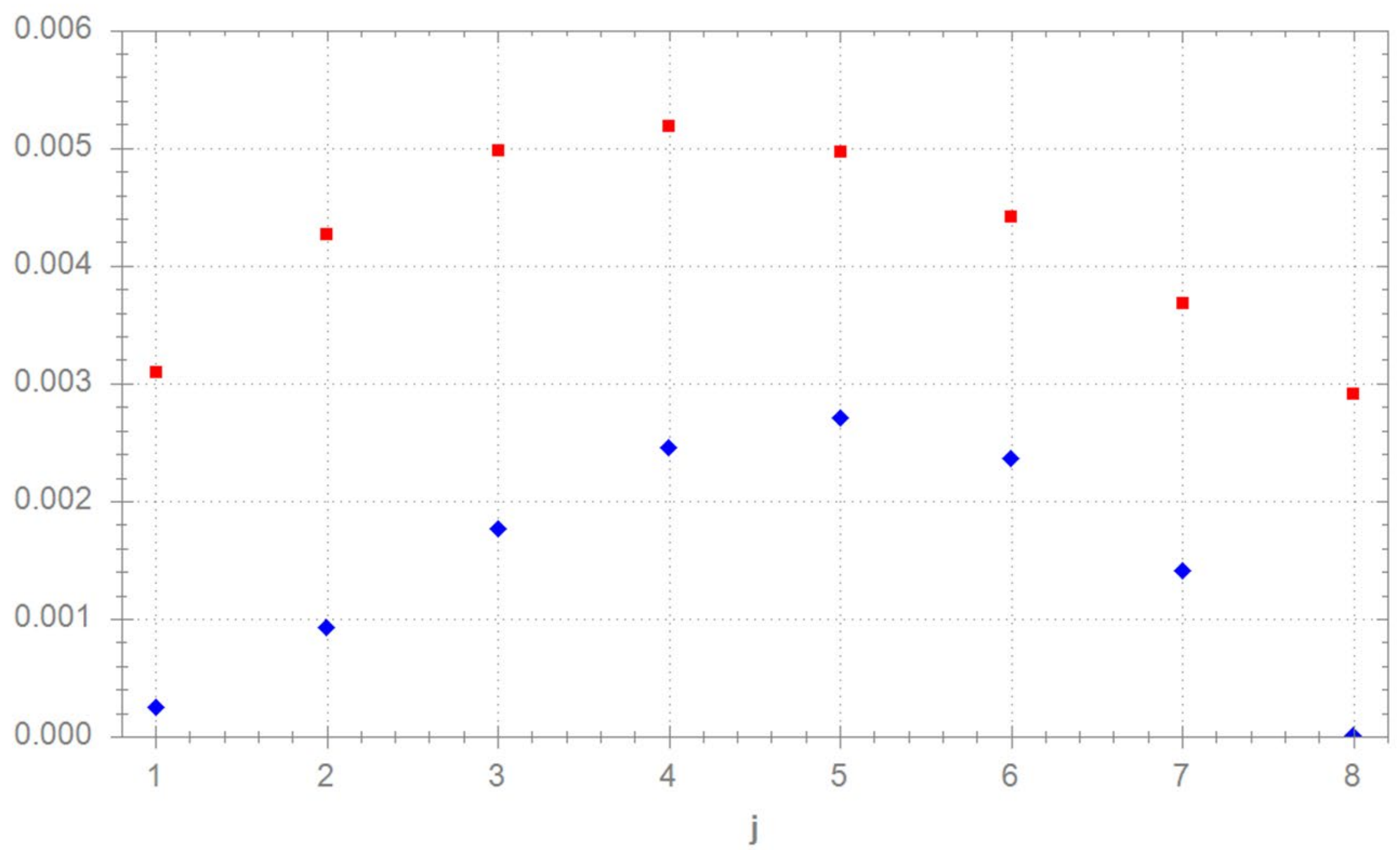

Figure 4.2

The differences between the left-hand and right-hand sides of (4.1) are strictly increasing in $j$.
$n = 1024,\ l = 256,\ u = 0,\ a = 8.$

Red boxes are $\dfrac{g(j+a+1)}{h(j+a+1)} - 1$.

Blue diamonds are $\left\{1 - \dfrac{g(j)}{h(j)}\right\}\left\{\dfrac{j(j+l+a+1)}{(j+l)(j+a+1)}\right\}\dfrac{\dbinom{2l-1+2j}{l-1+j}\dbinom{2l-1-2j}{l-1+j}}{\dbinom{2l+1+2a+2j}{l+a+j}\dbinom{2l-3-2a-2j}{l+a+j}}$.

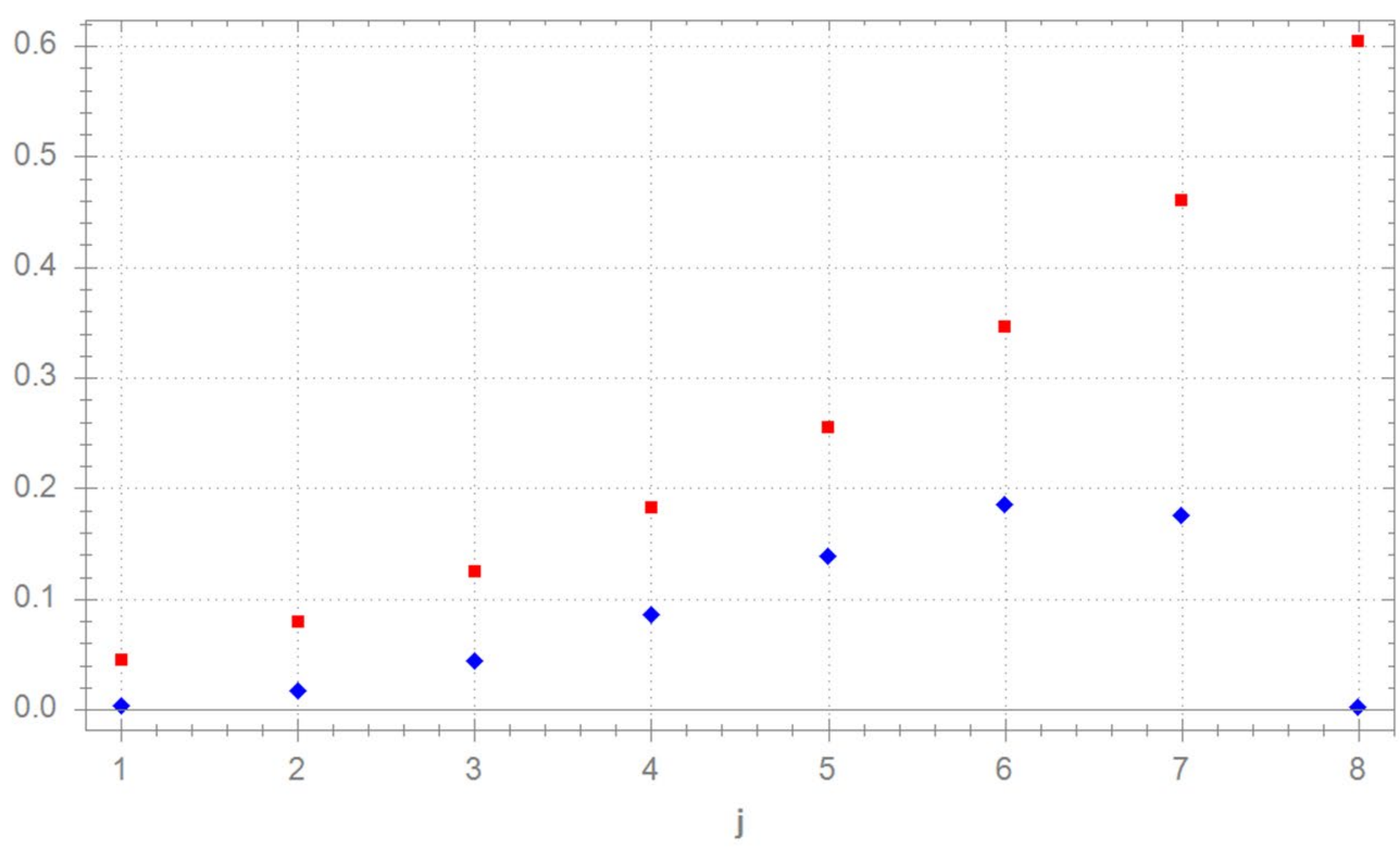

Figure 4.3

The differences between the left-hand and right-hand sides of (4.1) are strictly decreasing in $l$ $a=3,\ j=1,\ l=36,\ldots,63$.

Red boxes are $\dfrac{g(j+a+1)}{h(j+a+1)}-1$.

Blue diamonds are $\left\{1-\dfrac{g(j)}{h(j)}\right\}\left\{\dfrac{j(j+l+a+1)}{(j+l)(j+a+1)}\right\}\dfrac{\binom{2l-1+2j}{l-1+j}\binom{2l-1-2j}{l-1+j}}{\binom{2l+1+2a+2j}{l+a+j}\binom{2l-3-2a-2j}{l+a+j}}$.

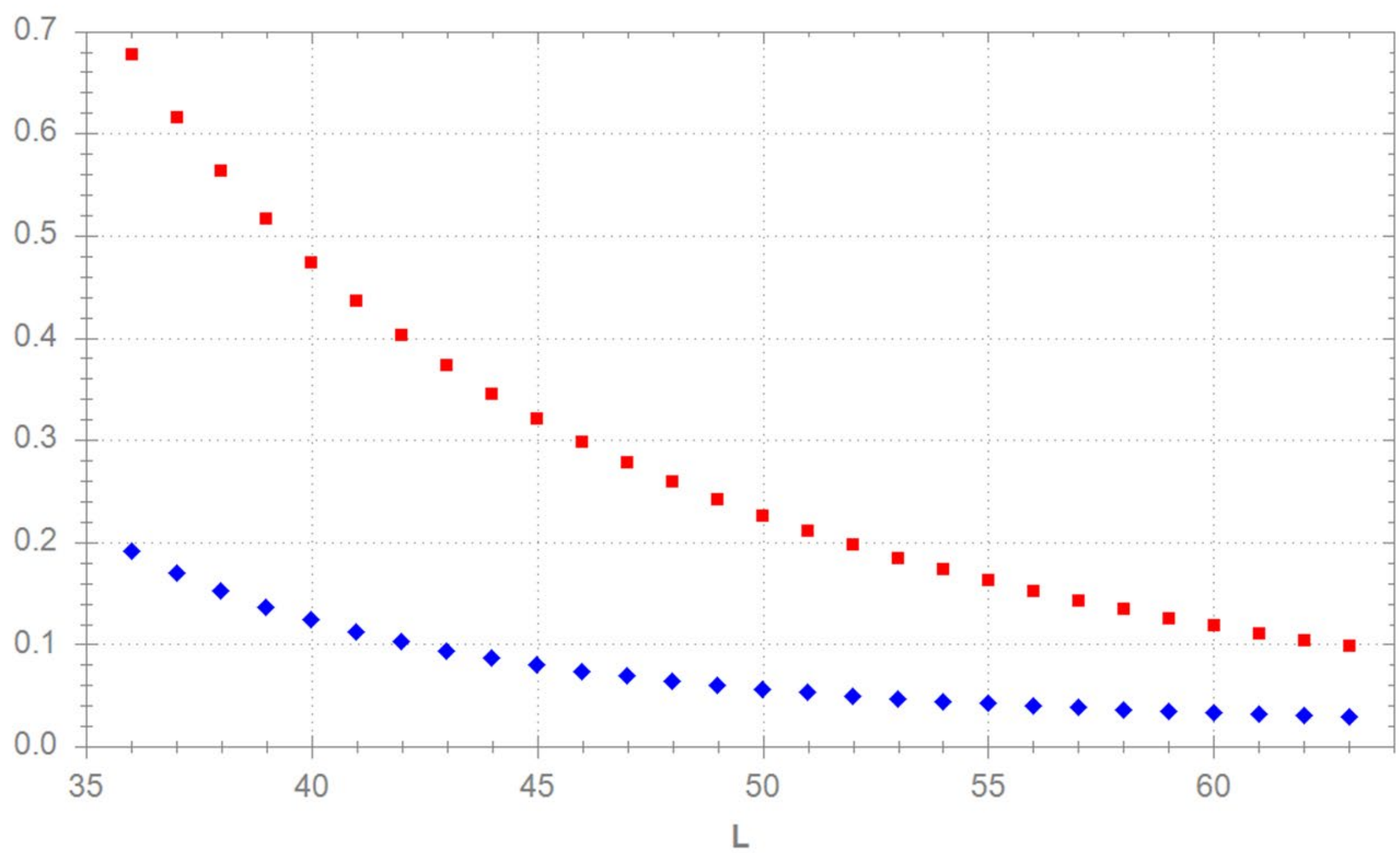

Figure 4.4

The differences between the left-hand and right-hand sides of (4.1) are strictly decreasing in $l$
$a=3,\ j=2,\ l=36,\ldots,63$.

Red boxes are $\dfrac{g(j+a+1)}{h(j+a+1)}-1$.

Blue diamonds are $\left\{1-\dfrac{g(j)}{h(j)}\right\}\left\{\dfrac{j(j+l+a+1)}{(j+l)(j+a+1)}\right\}\dfrac{\dbinom{2l-1+2j}{l-1+j}\dbinom{2l-1-2j}{l-1+j}}{\dbinom{2l+1+2a+2j}{l+a+j}\dbinom{2l-3-2a-2j}{l+a+j}}$.

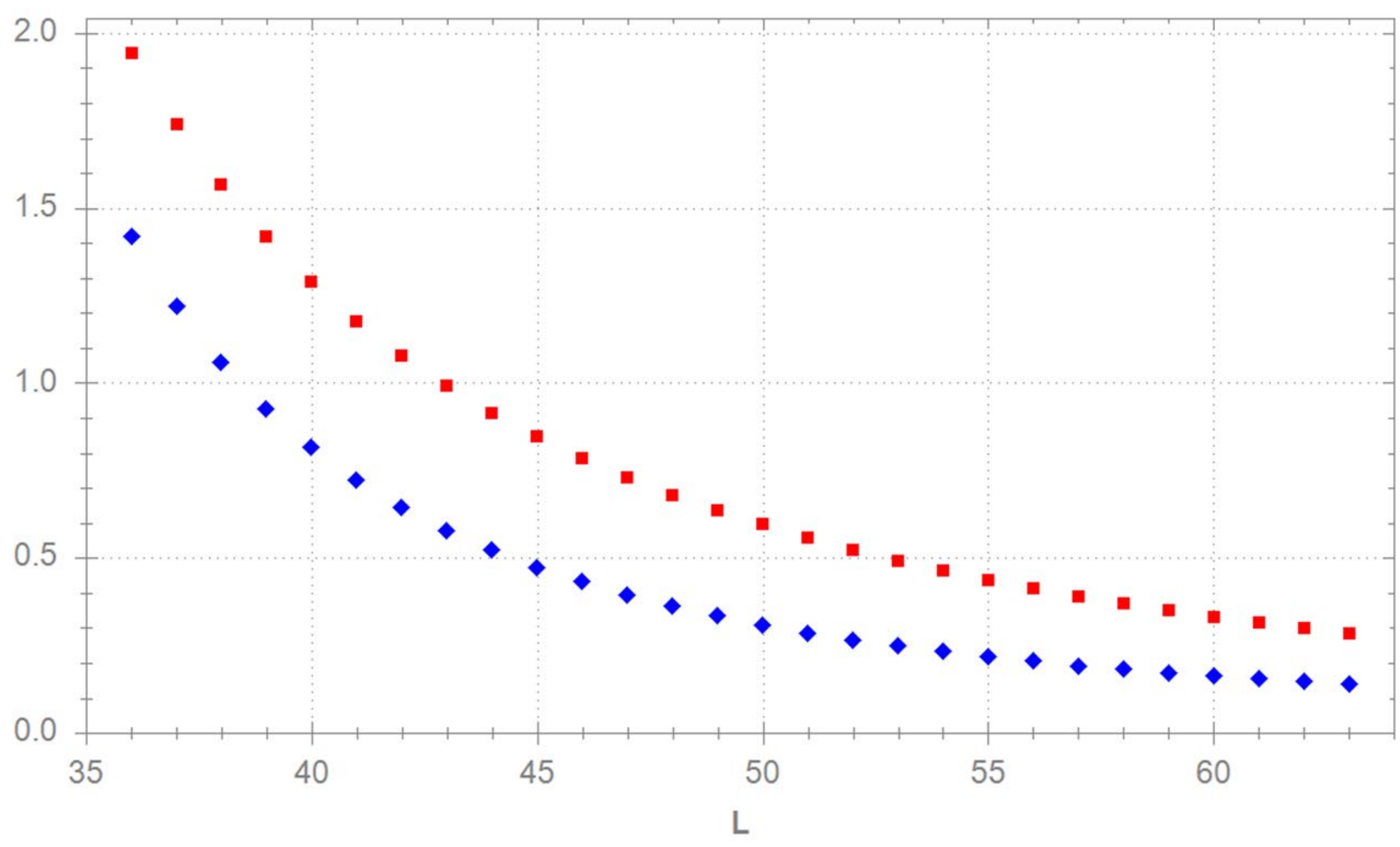

Figure 4.5

The differences between the left-hand and right-hand sides of (4.1) are strictly decreasing in $l$
$a=3,\ j=3,\ l=36,\ldots,63$.

Red boxes are $\dfrac{g(j+a+1)}{h(j+a+1)}-1$.

Blue diamonds are $\left\{1-\dfrac{g(j)}{h(j)}\right\}\left\{\dfrac{j(j+l+a+1)}{(j+l)(j+a+1)}\right\}\dfrac{\dbinom{2l-1+2j}{l-1+j}\dbinom{2l-1-2j}{l-1+j}}{\dbinom{2l+1+2a+2j}{l+a+j}\dbinom{2l-3-2a-2j}{l+a+j}}$.

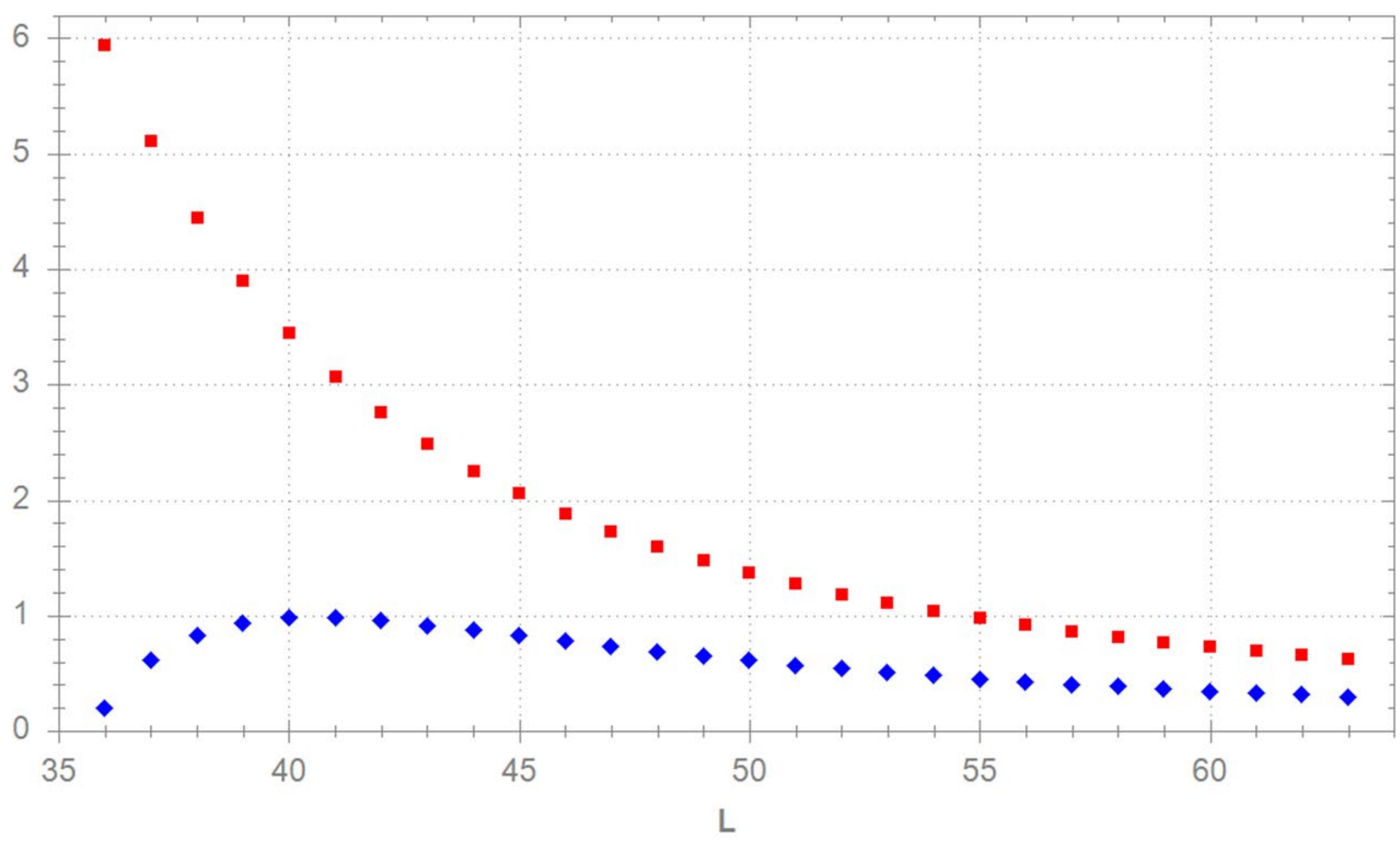

Figure 4.6

Differences $\Delta(a, j, l)$ between the left-hand and right-hand sides of (4.1) for (original) $a=1,\ldots,10$ in the most stringent case, $j=1$, $l=4(a+1)^2-1$).

Red boxes are $\Delta(a, 1, 4(a+1)^2-1)$.

Blue diamonds are $\{1/(a+1)^2\} + (\frac{19}{12} - \frac{e}{2})/(a+1)^3$.

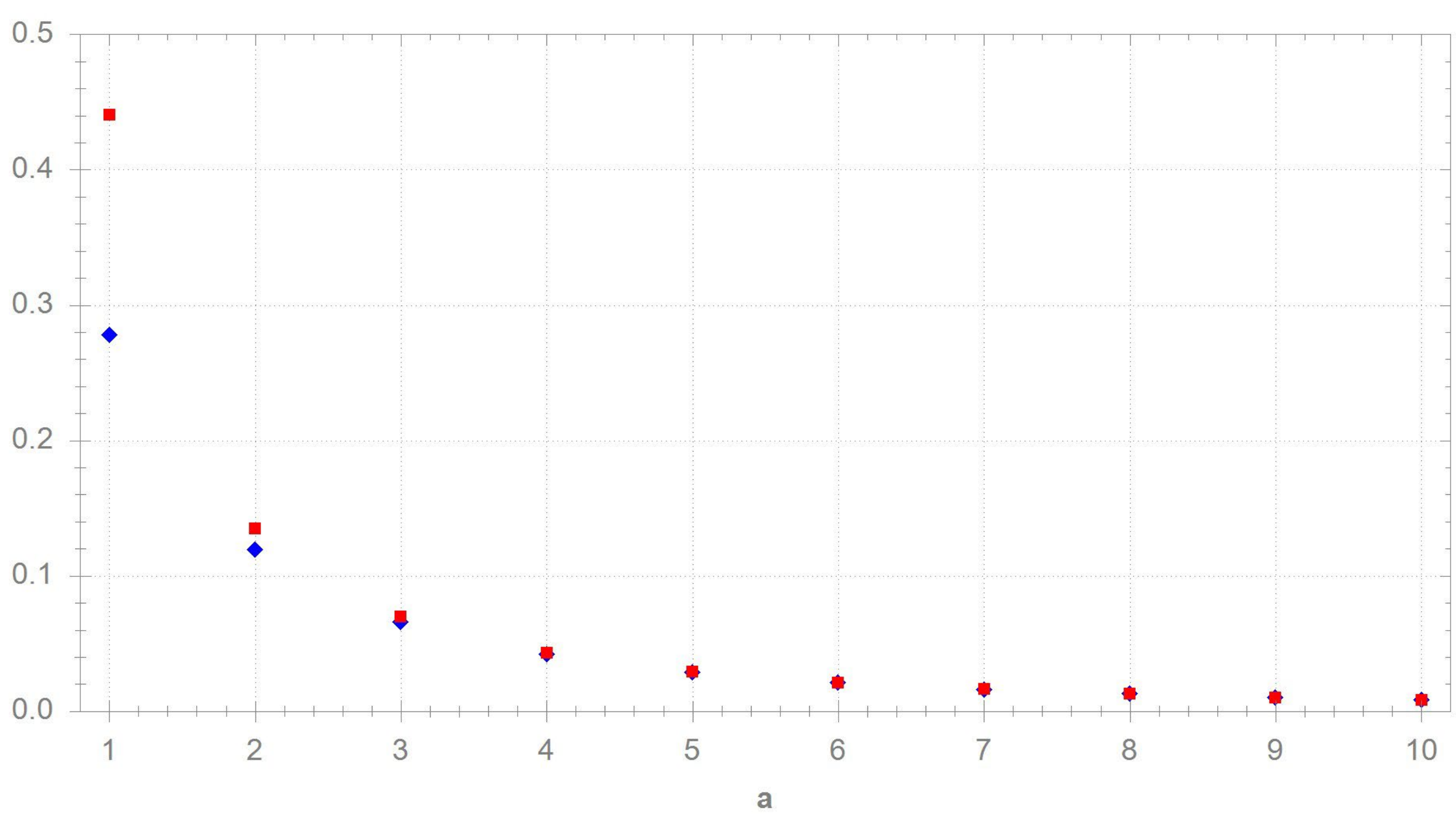

Figure 4.7a

Plot of $\Delta(a, j, l)$ versus $j$ and $l$ in the case $a=5$.

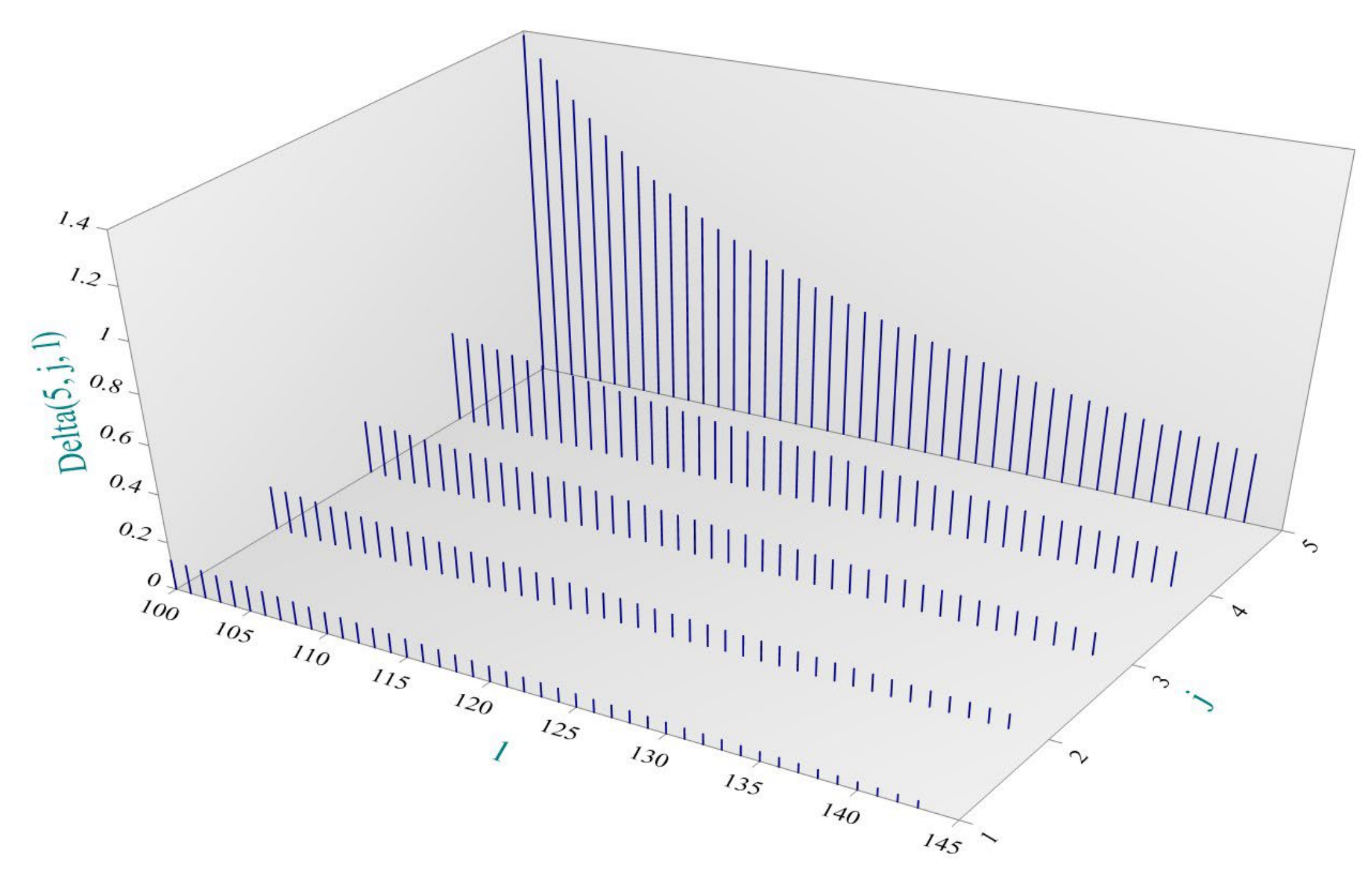

Figure 4.7b

Plot of $\Delta(a, j, l)$ versus $j$ and $l$ in the case $a=8$.

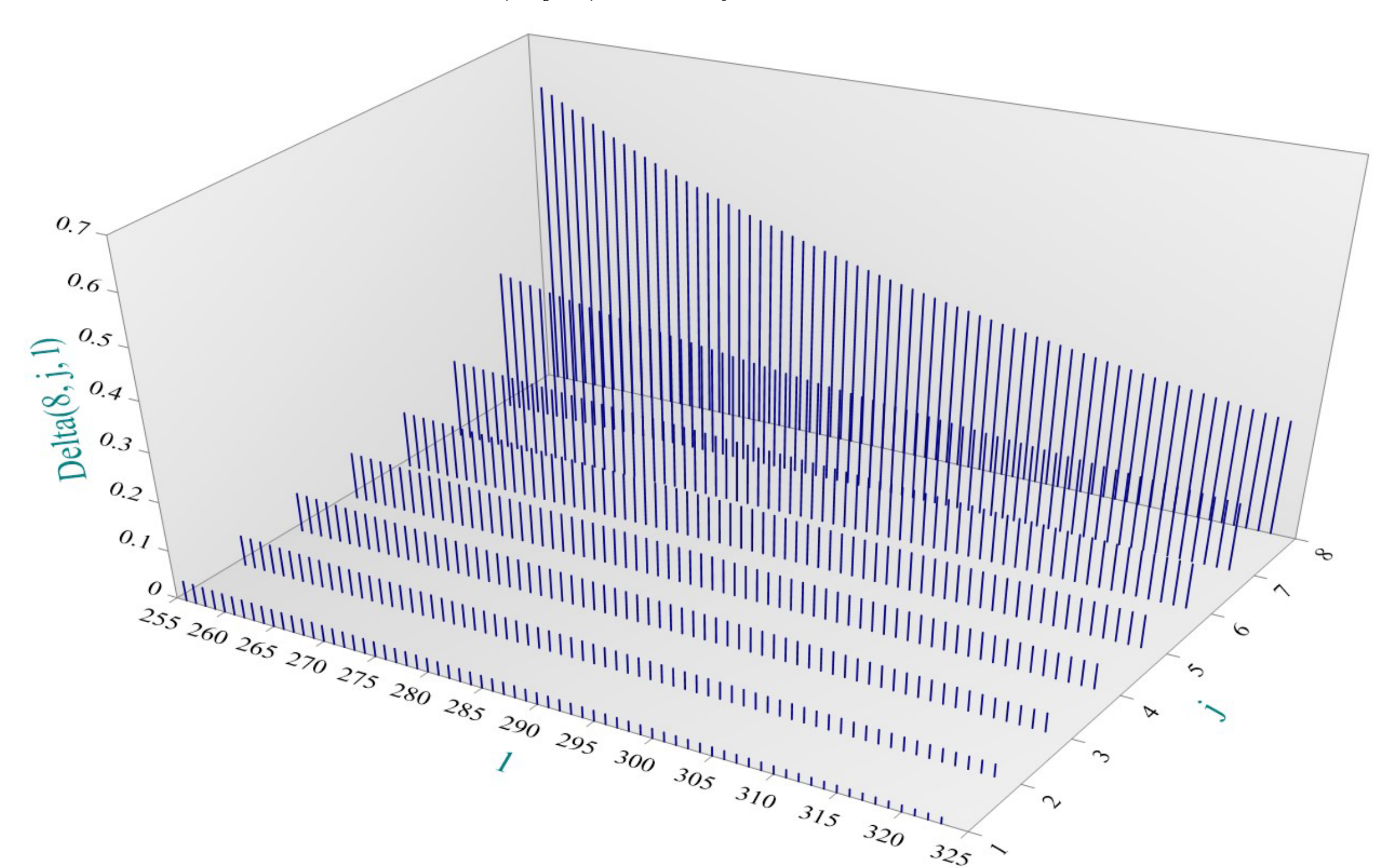